\documentclass[12pt, reqno, twoside]{amsart}
\usepackage{amsmath, amsthm, amscd, amsfonts, amssymb, graphicx, color}
\usepackage[bookmarksnumbered, colorlinks, plainpages]{hyperref}
\hypersetup{colorlinks=true,linkcolor=red, anchorcolor=green, citecolor=cyan, urlcolor=red, filecolor=magenta, pdftoolbar=true}
\usepackage{enumerate,amsfonts,mathtools,esint}
\usepackage{mathrsfs}
\usepackage{fancyhdr}
\usepackage{enumitem}
\usepackage{array}
\usepackage{hyperref}

\usepackage[toc,page]{appendix}
\usepackage{tabularx}
\usepackage{graphicx}

\emergencystretch=\maxdimen
\usepackage{float}

\usepackage{tikz}

\usepackage{bm}

\DeclareFontFamily{OT1}{rsfs}{}
\DeclareFontShape{OT1}{rsfs}{m}{n}{ <-7> rsfs5 <7-10> rsfs7 <10-> rsfs10}{}
\DeclareMathAlphabet{\mycal}{OT1}{rsfs}{m}{n}

\def\b0{{\mathbf{0}}}

\newtheorem{theorem}{Theorem}[section]

\newtheorem{remark}{Remark}[section]
\newtheorem{lemma}{Lemma}[section]

\numberwithin{equation}{section}

\makeatletter
\def\hlinewd#1{%
\noalign{\ifnum0=`}\fi\hrule \@height #1 %
\futurelet\reserved@a\@xhline}
\makeatother
\newcolumntype{?}[1]{!{\vrule width #1}}

\newcounter{marnote}

\tikzset{
	MyPersp/.style={scale=1.8,x={(-0.8cm,-0.4cm)},y={(0.8cm,-0.4cm)},
    z={(0cm,1cm)}},
	MyPoints/.style={fill=white,draw=black,thick}
	}

\begin{document}

\setcounter{page}{1}

\title[Curvature conditions]{Curvature conditions, Liouville-type theorems and Harnack inequalities for a nonlinear 
parabolic equation on smooth metric measure spaces}

\author[A. Taheri]{Ali Taheri}
\author[V. Vahidifar]{Vahideh Vahidifar}
\address{School of Mathematical and  Physical Sciences, 
University of Sussex, Falmer, Brighton, United Kingdom.}
\email{\textcolor[rgb]{0.00,0.00,0.84}{a.taheri@sussex.ac.uk}} 
\address{School of Mathematical and  Physical Sciences, 
University of Sussex, Falmer, Brighton, United Kingdom.}
\email{\textcolor[rgb]{0.00,0.00,0.84}{v.vahidifar@sussex.ac.uk}}


\begin{abstract}
In this paper we prove gradient estimates of both elliptic and parabolic types, specifically, of Souplet-Zhang, 
Hamilton and Li-Yau types for positive smooth solutions to a class of nonlinear parabolic equations involving 
the Witten or drifting Laplacian on smooth metric measure spaces. These estimates 
are established under various curvature conditions and lower bounds on the generalised Bakry-\'Emery Ricci 
tensor and find utility in proving elliptic and parabolic Harnack-type inequalities as well as general Liouville-type and 
other global constancy results. Several applications and consequences are presented and discussed. 
\end{abstract}

\maketitle
 
{ 
\hypersetup{linkcolor=black}
\tableofcontents
}

\section{Introduction} \label{sec1}

Unravelling the interplay between the analytic, stochastic and geometric properties of Riemannian manifolds on the one hand and 
the dynamics of the solutions to nonlinear evolution equations on the other lies at the heart of geometric analysis. Some particular 
questions include: short and long time existence, smoothing properties, convergence to equilibrium, spectral gaps, 
blow-ups, self-similar solutions and their profiles, analysis of singularities, critical exponents, 
description and classification of ancient and eternal solutions as well as Liouville or global constancy type 
results to mention a few.

The ubiquitous gradient estimates play a significant role here and their derivation is often the main step in 
confronting such challenging tasks. For parabolic and heat-type equations, 
there are several types of gradient estimates and techniques in literature, each having their 
own strength and utility. In this paper we extend such techniques and results to the context of nonlinear equations 
on smooth metric measure spaces where understanding the interaction between the 
nonlinearity and curvature on the one hand and the different terms in the estimate on the other is the main driving force. 
({\it See} \cite{[Ha93],[LY86]}, \cite{AM, Bid, Giaq,GidSp, SZ, Taheri-GE-1, Taheri-GE-2, [Wang], [Wu15], [Wu18]},   
\cite{Aub, [Ba1], BDM, Gr, [LiP],SchYau,Taheri-book-two, [QZh]} for more).

Towards this end let $(M,g)$ be a smooth Riemannian manifold of dimension $n \ge 2$ and let 
$d\mu=e^{-\phi} dv_g$ denote the weighted measure associated with a potential 
$\phi$ and the Riemannian volume measure $dv_g$ on $M$. 
The triple $(M, g, d\mu)$ is referred to as a smooth metric space or a weighted manifold. 
In this paper we develop gradient estimates of elliptic and parabolic types, more specifically, 
of Souplet-Zhang, Hamilton and Li-Yau types, along with some of their implications, for positive 
smooth solutions $w=w(x,t)$ to the nonlinear parabolic equation on $(M,g,d\mu)$ given by:
\begin{align} \label{eq11}
({\mathscr P}) \qquad \qquad \frac{\partial w}{\partial t} = \Delta_\phi w 
+ {\mathscr G}(t, x, w), \qquad (x,t) \in M \times (-\infty, \infty).  
\end{align}

The operator $\Delta_\phi$ appearing on the right in \eqref{eq11} is the Witten Laplacian associated with the smooth 
metric measure space $(M,g,d\mu)$ (also called the weighted or drifting or $\phi$-Laplacian) whose action on functions 
$v \in \mathscr{C}^2(M)$ is given by 
\begin{equation} \label{Lf definition}
\Delta_\phi v = e^\phi {\rm div} (e^{-\phi} \nabla v)= \Delta v - \langle \nabla \phi, \nabla v \rangle.   
\end{equation}

Here $\Delta, {\rm div}$ and $\nabla$ are the usual Laplace-Beltrami, divergence and gradient operators 
associated with the metric $g$. The nonlinearity ${\mathscr G}={\mathscr G}(t,x,w)$ is a sufficiently 
smooth function of the space-time variables $(t,x)$ as well as the dependent 
variable $w$. The general form of the assumed nonlinearity here enables us to more clearly 
examine the form and extent in which it influences the various gradient estimates and subsequently 
the implication it has on the qualitative properties of solutions.

Gradient estimates for positive solutions to linear 
and nonlinear heat type equations have been studied extensively starting from the seminal paper of Li and Yau \cite{[LY86]} 
({\it see} also \cite{[LiP]}). In the nonlinear setting perhaps the first equation to be considered is the one with a logarithmic 
type nonlinearity ({\it see}, e.g., \cite{LiJ91}) [specifically, ${\mathscr G} = {\mathsf A}(x,t) w \log w$ in $({\mathscr P})$] given by
\begin{equation} \label{eq1.4}
\frac{\partial w}{\partial t} = \Delta_\phi w + {\mathsf A}(x,t) w \log w. 
\end{equation}
The interest in such problems originates partly from its natural links with gradient Ricci solitons and partly from links with 
geometric and functional inequalities on manifolds, notably, the logarithmic Sobolev and energy-entropy 
inequalities \cite{Bak, Gross, Taheri-book-one, Taheri-book-two, VC, Wang}. Recall that a Riemannian manifold $(M,g)$ is said to be a gradient 
Ricci soliton if there there exists a smooth function $\phi$ on $M$ and a constant $\lambda \in \mathbb{R}$ 
such that ({\it cf.} \cite{Cao, Lott})
\begin{equation}\label{eq1.5}
{\mathscr Ric}_\phi(g) = {\mathscr Ric}(g) + \nabla \nabla \phi = \lambda g.
\end{equation}
The notion is a generalisation of an Einstein manifold and has a fundamental role in the analysis of singularities 
of the Ricci flow \cite{CaoChow, Chow, Ham, [QZh]}. Taking trace from both sides of \eqref{eq1.5} and using the contracted 
Bianchi identity leads one to a simple form of \eqref{eq1.4} with constant coefficients.

Another prominent class of nonlinear equations rooted in conformal geometry and studied extensively in this setting are 
Yamabe-type equations ({\it see}, \cite{Bid, GidSp, Lee, Mast}). In the context of smooth metric measure spaces these 
equations can be broadly formulated as ({\it see}, e.g., \cite{Case, Dung, [Wu18]}) 
\begin{align} \label{eq11s}
\frac{\partial w}{\partial t} = \Delta_\phi w + {\mathsf A}(x,t) w^p + {\mathsf B}(x,t) w. 
\end{align}
Incidentally, the case ${\mathsf A} \equiv -1$, ${\mathsf B} \equiv 1$, $p=3$ [${\mathscr G}(w) = w-w^3$]  
is the Allen-Cahn equation and the case ${\mathsf A} \equiv -c$, 
${\mathsf B} \equiv c$, $p=2$ [${\mathscr G}(w)=cw(1-w)$ with $c>0$] is the Fisher-KPP equation ({\it cf.} \cite{AllC, Fish, KPP}). 
Both these equations have been studied extensively in recent years due to the significance of the phenomenon they model and 
their huge applications in physics and other sciences. (For various geometric estimates and their consequences 
{\it see} \cite{BM, CLPW} and the references therein). 
A far reaching generalisation of \eqref{eq11s} with a superposition of power-like nonlinearities consist of equations in the form 
\begin{align}
\frac{\partial w}{\partial t} = \Delta_\phi w + \sum_{j=1}^d {\mathsf A}_j(x,t) w^{p_j} + \sum_{j=1}^d {\mathsf B}_j(x,t) w^{q_j}. 
\end{align}
Here ${\mathsf A}_j$, ${\mathsf B}_j$ (with $1 \le j \le d$) are sufficiently smooth space-time dependent coefficients and 
$p_j \ge 0$, $q_j \le 0$  real exponents ({\it see} \cite{Taheri-GE-1, Taheri-GE-2}). 
Other classes of equations generalising the above and close to \eqref{eq1.4} and \eqref{eq11s} appear in the form 
({\it see}, e.g., \cite{CGS, GKS, Taheri-GE-1, Taheri-GE-2, [Wu15]}) 
\begin{equation}
\frac{\partial w}{\partial t} = \Delta_\phi w + {\mathsf A}(x,t) \Gamma(\log w) w^p + {\mathsf B}(x,t) w^q + {\mathsf C}(x,t) w, 
\end{equation}
with $p$, $q$ real exponents, ${\mathsf A}$, ${\mathsf B}$, ${\mathsf C}$ sufficiently smooth space-time dependent coefficients and 
$\Gamma \in {\mathscr C}^1(\mathbb{R}, \mathbb{R})$. Some cases of particular interest for $\Gamma=\Gamma(s)$ include power function, 
for example, $s^\alpha$ (integer $\alpha \ge 1$), $|s|^\alpha$ or $|s|^{\alpha-1} s$ (real $\alpha >1$) with different 
sign-changing, growth and singular behaviour as $s \to \pm \infty$, or a superposition of such nonlinearities 
\cite{Taheri-GE-2}. 

Another case of recent interest arises from iterated logarithms ({\it cf.}  \cite{CGS}) associated with a string of parameters 
$d$, $k_1, \cdots, k_d \in {\mathbb N}$ and $\beta_1, \cdots, \beta_d \in {\mathbb R}$, specifically, 
\begin{equation}
\Gamma_{k_1, \dots, k_d}^{\beta_1, \dots, \beta_d} (\log w) = |\log_{k_1} w|^{\beta_1} |\log_{k_2} w|^{\beta_2} \dots |\log_{k_d} w|^{\beta_d}, 
\end{equation}
where $\log_k w = \log \log_{k-1} w$ for $k \ge 2$ and $\log_1 w = \log w$. This can be considered but with due care, e.g., 
subject to the assumption of $w$ being sufficiently large, specifically, with respect to iterated exponentials of 
$k_1, \dots, k_d$ (as otherwise the repeated logarithm is meaningless due to the 
possibility of $\log_{k-1} w$ being non-positive hence making $\log_k w$ undefined). 
Naturally one can also consider 
variations of the same theme, for example, by replacing $\log_k$ with either of its close relatives
\begin{align}
&\log^{||}_k w = \log^{||} \log^{||}_{k-1} w \mbox{ for $k \ge 2$}, \qquad \,\, \log^{||}_1 w = |\log w|, \\ 
&\log^+_k w = \log^+ \log^+_{k-1} w \mbox{ for $k \ge 2$}, \qquad \log^+_1 w = 1+ [\log w]_+.
\end{align}
However, one needs to observe that the function $\Gamma$ thus obtained is only ${\mathscr C}^1$ 
outside a discrete set (the zero sets of the functions $\log_{k-1} w$) and hence does not lie in 
${\mathscr C}^1({\mathbb R}, {\mathbb R})$ as required.

Another related and yet more general form of Yamabe-type equations is the Einstein-scalar field Lichnerowicz 
equation ({\it see} Choquet-Bruhat \cite{CBY}, Chow \cite{Chow} and Zhang \cite{[QZh]}). 
In the context of smooth metric measure spaces a generalisation of the Einstein-scalar field Lichnerowicz 
equation with space-time dependent coefficients can be described as:  
\begin{align}
\frac{\partial w}{\partial t} &= \Delta_\phi w + {\mathsf A}(x,t) w^p + {\mathsf B}(x,t) w^q + {\mathsf C}(x,t) w \log w, \\
\frac{\partial w}{\partial t} &= \Delta_ \phi w + {\mathsf A}(x,t) e^{2w} + {\mathsf B}(x,t) e^{-2w} + {\mathsf C}(x,t). 
\end{align}
For gradient estimates, Harnack inequalities, Liouville type theorems and other related results in this direction 
see \cite{Dung, Ma, Song, Taheri-GE-1, Taheri-GE-2, TVahNonlin,TVahNA, [Wu18]} and the references therein.

\qquad \\
{\bf Plan and layout of paper.} Let us conclude this introduction by briefly describing the plan of the paper. 
In Section \ref{sec1-1} we present some background material on smooth metric measure spaces and the 
associated generalised Ricci curvature tensors as required for the development of the paper. 
In Section \ref{sec2} we present the main results of the paper along with some related discussion. 
The remainder of the paper is then devoted to the detailed proofs of these results. In Section \ref{sec3} we 
present the proof of the local Souplet-Zhang type gradient estimate in Theorem \ref{thm1} and in Section \ref{sec4} 
we give the proof of the local and global elliptic Harnack inequality in Theorem \ref{cor Harnack}. In Section \ref{sec5} we 
present the proof of the local Hamilton-type  gradient estimate in Theorem \ref{thm18} and in Section \ref{sec6} 
we give the proof of the various Liouville-type results formulated in Theorems \ref{Liouville one thm} and  
\ref{Liouville two thm}. In Section \ref{sec7} we present the proof of the local 
Li-Yau differential Harnack estimate in Theorem \ref{thm7.1} and in Section \ref{sec8} we establish the local and global 
parabolic Harnack inequalities as formulated in Theorem \ref{thm38}. 
Section \ref{sec9} is devoted to the proof of the general Liouville result in Theorem 
\ref{coroLiouville} and its subsequent corollaries.

For the convenience  of the reader and future reference we describe at the end of this section some of the key notation and quantities 
used in the paper.

\qquad \\
{\bf Notation.} 
We denote by $d=d(x,y)$ the geodesic distance between $x$, $y$ in $M$. 
For the sake of the local estimates below we typically fix a reference point $x_0$ in $M$ and denote by $r=r(x)$ the geodesic radial 
variable measuring the distance between $x$, $x_0$ in $M$. We write ${\mathscr B}_R(x_0)$ 
for the closed geodesic ball in $M$ centred at $x_0$ with radius $R>0$ and we write 
$Q_{R,T}(x_0)$ for the compact parabolic space-time cylinder with lower base 
${\mathscr B}_R(x_0) \times \{t_0-T\}$ for $t_0 \in {\mathbb R}$ and height $T>0$, specifically, 
\begin{equation} \label{parabolic-cylinder-Q}
Q_{R,T}(x_0)={\mathscr B}_{R}(x_0) \times [t_0-T, t_0] \subset M \times (-\infty, \infty).
\end{equation} 
When $t_0=T>0$ it is more convenient to write \eqref{parabolic-cylinder-Q} as
\begin{equation}
H_{R,T}(x_0) = {\mathscr B}_R(x_0) \times [0, T] \subset M \times [0, T].
\end{equation} 
When the choices of $x_0$ and $t_0$ are clear from the context, 
in the interest of brevity we simply write ${\mathscr B}_R$, $Q_{R, T}$ and $H_{R,T}$. 
We write $X_+=\max(X, 0)$ and $X_-=\max(-X, 0)$.

\section{Preliminaries and Background} \label{sec1-1}

A smooth metric measure space or a weighted manifolds is a triple $(M, g, d\mu)$ 
where $(M,g)$ is a smooth Riemannian manifold of dimension $n \ge 2$ and 
$d\mu=e^{-\phi} dv_g$ is the weighted measure associated with the potential 
$\phi$ (that is, the positive weight $\omega=e^{-\phi}$) and the Riemannian volume 
measure $dv_g$ on $M$. The notion is an important one in geometric analysis 
and arises in various contexts ({\it see}, e.g., \cite{Gr, [LLi4], [LLi1], [LLi2], [Li0], [Pe02], [QZh]}).

Associated with the triple $(M, g, d\mu)$ there is a $\phi$-Laplacian 
(also called the weighted, drifting or Witten Laplacian) defined by \eqref{Lf definition}. 
The $\phi$-Laplacian is a symmetric diffusion operator with respect to the invariant weighted measure 
$d\mu$ and reduces to the ordinary Laplacian (the Laplace-Beltrami operator) precisely when the potential $\phi$ 
is a constant function. By an application of the integration by parts formula it can be seen that for any compactly 
supported smooth functions  $u, w \in {\mathscr C}_0^\infty(M)$ we have  
\begin{equation}
\int_M e^{-\phi}  w \Delta_\phi u \, dv_g 
= - \int_M e^{-\phi} \langle \nabla u, \nabla w \rangle \, dv_g 
= \int_M e^{-\phi} u \Delta_\phi w \, dv_g.
\end{equation}

As for the geometry and curvature properties of the triple $(M,g,d\mu)$ we introduce the generalised Ricci curvature 
tensor field on $M$ by writing    
\begin{equation} \label{Ricci-m-f-intro}
{\mathscr Ric}^m_\phi(g) = {\mathscr Ric}(g) + \nabla\nabla \phi - \frac{\nabla \phi \otimes \nabla \phi}{m-n}, 
\end{equation}
(the Bakry-\`Emery $m$-Ricci curvature) where ${\mathscr Ric}(g)$ is the standard Riemannain Ricci curvature of $g$, 
$\nabla \nabla \phi={\rm Hess}(\phi)$ denotes the Hessian of $\phi$, and $m \ge n$ is a constant ({\it see}, e.g., \cite{BD, BE, Bak, Lott, VC, LFWang-a, LFWang-b}). 
For the sake of clarity we point out that when 
$m=n$, by convention, $\phi$ is only allowed to be a constant, thus giving ${\mathscr Ric}^m_\phi(g)={\mathscr Ric}(g)$, 
whereas, we also allow for $m = \infty$ in which case by formally passing to the limit in \eqref{Ricci-m-f-intro} we get, 
\begin{equation} \label{Ricf def eq} 
{\mathscr Ric}_\phi(g) = {\mathscr Ric}(g) + \nabla\nabla \phi := {\mathscr Ric}^\infty_\phi(g).
\end{equation}

The weighted Bochner-Weitzenb\"ock formula in this context asserts that for every function $w \in {\mathscr C}^3(M)$ we have the pointwise differential identity
\begin{equation} \label{Bochner}
\frac{1}{2} \Delta_\phi |\nabla w|^2 
= |\nabla\nabla w|^2 + \langle \nabla w, \nabla \Delta_\phi w \rangle 
+ {\mathscr Ric}_\phi (\nabla w, \nabla w).
\end{equation}

Making note of the basic inequality $(\Delta w)^2/n \le |\nabla \nabla w|^2$ and by recalling 
the defining relation $\Delta_\phi w = \Delta w - \langle \nabla \phi, \nabla w \rangle$ it is then easily seen that 
\begin {align} \label{Eq-1.7z}
|\nabla\nabla w|^2 + \frac{\nabla \phi \otimes \nabla \phi}{m-n} (\nabla w, \nabla w)
&= |\nabla\nabla w|^2 + \frac{\langle \nabla \phi , \nabla w \rangle^2}{m-n} \nonumber \\
&\ge \frac{(\Delta w)^2}{n} + \frac{\langle \nabla \phi, \nabla w \rangle ^2}{m-n} \nonumber \\
&\ge \frac{(\Delta w - \langle \nabla \phi, \nabla w \rangle)^2}{m} = \frac{(\Delta_\phi w)^2}{m}.  
\end{align}
Subsequently it follows from \eqref{Ricci-m-f-intro}, \eqref{Bochner} and \eqref{Eq-1.7z} that here we have the inequality 
\begin{equation} \label{Bochner-m-inequality}
\frac{1}{2} \Delta_\phi |\nabla w|^2 
\ge \frac{1}{m} (\Delta_\phi w)^2 + \langle \nabla w, \nabla \Delta_\phi w \rangle 
+ {\mathscr Ric}^m_\phi (\nabla w, \nabla w).
\end{equation}

Now a curvature lower bound in the form ${\mathscr Ric}^m_\phi(g) \ge k g$ (with $k \in {\mathbb R}$) 
implies that the diffusion operator $L=\Delta_\phi$ satisfies the curvature-dimension 
condition ${\rm CD}(k,m)$ ({\it see} \cite{BE, Bak, Lott, VC, Wang}). To elaborate on this last point 
further, we first recall the definition of the carre du champ operator $\Gamma[L]$ 
associated with a Markov diffusion operator $L$. This is given explicitly by the formulation
\begin{equation} \label{Gamma-One-L-eq}
\Gamma[L] (u,v) = \Gamma_1[L] (u,v) := \frac{1}{2} [L(uv) - u L v - v Lu]. 
\end{equation} 

Higher order iterates are then defined inductively by replacing the product operation with $\Gamma$ respectively 
({\it cf.} \cite{Bak}). As a matter of fact, the second order iterated carre du champ operator $\Gamma_2[L]$, can 
be seen to be given by, 
\begin{equation} \label{Gamma-Two-L-eq}
\Gamma_2[L] (u,v) := \frac{1}{2} [ L \Gamma(u,v) - \Gamma (u, Lv) - \Gamma (Lu, v)]. 
\end{equation}

By a direct calculation, it is now observed that for the diffusion operator $L=\Delta_\phi$ as in \eqref{Lf definition} and with 
$\Gamma[L]$ and $\Gamma_2[L]$ as in \eqref{Gamma-One-L-eq} and \eqref {Gamma-Two-L-eq}, we have the first 
and second order relations 
\begin{equation} \label{Un-iterated-CDC}
\Gamma[L] (w,w) = |\nabla w|^2,
\end{equation}
and 
\begin{equation} \label{iterated-CDC}
\Gamma_2 [L] (w,w) = \frac{1}{2} L |\nabla w|^2 - \langle \nabla w, \nabla L w \rangle, 
\end{equation} 
respectively. Hence, in light of \eqref{Bochner-m-inequality}, \eqref{Un-iterated-CDC} and \eqref{iterated-CDC} it follows from the 
curvature lower bound ${\mathscr Ric}^m_\phi(g) \ge k g$ that here the 
iterated carre du champ operator $\Gamma_2[L]$ satisfies the inequality 
\begin{align} \label{CD-mk}
\Gamma_2[L](w,w) &=  \frac{1}{2} L |\nabla w|^2 - \langle \nabla w, \nabla L w \rangle \nonumber \\
&\ge \frac{1}{m} (L w)^2 + {\mathscr Ric}^m_\phi (\nabla w, \nabla w) \nonumber \\
&\ge \frac{1}{m} (L w)^2 + k |\nabla w|^2 = \frac{1}{m} (Lw)^2 + k \Gamma[L](w,w). 
\end{align}
In contrast, a curvature lower bound in the form ${\mathscr Ric}_\phi(g) \ge k g$ implies the curvature-dimension condition 
${\rm CD}(k,\infty)$, in the sense that by virtue of \eqref{Bochner}, \eqref{Un-iterated-CDC} and \eqref{iterated-CDC} one only has
\begin{align} \label{CD-inftyk}
\Gamma_2[L](w,w) &=  \frac{1}{2} L |\nabla w|^2 - \langle \nabla w, \nabla L w \rangle \nonumber \\
&= |\nabla \nabla w|^2 + {\mathscr Ric}_\phi (\nabla w, \nabla w) 
\ge k |\nabla w|^2 = k \Gamma[L](w,w). 
\end{align}

\section{Statement of the main results} \label{sec2}

Here we present the main results of the paper leaving the proofs, technicalities and further details involved 
to the subsequent sections. For the convenience of the reader and clarity of presentation, below we have 
grouped the results into four subsections.

\subsection{Two gradient estimates of elliptic types for $({\mathscr P})$}

We begin by presenting two estimates of elliptic type for positive solutions to $({\mathscr P})$. 
These are estimates of Souplet-Zhang and Hamilton types respectively and the task is to exploit 
the role of the nonlinearity and the geometry of the triple $(M,g, d\mu)$. Here we fix a reference 
point $x_0 \in M$ and $t_0 \in {\mathbb R}$, $R>2$, 
$T>0$. Since the solution $w=w(x,t)>0$ and the parabolic space-time cylinder 
$Q_{R,T}(x_0) \subset M \times (-\infty, \infty)$ is compact, $w$ is bounded 
away from zero and bounded from above in $Q_{R,T}(x_0)$ and so the 
quantities involving the nonlinearity and its derivatives in \eqref{eq13} are finite.

\begin{theorem} \label{thm1}
Let $(M, g, d\mu)$ with $d\mu=e^{-\phi} dv_g$ be a complete smooth metric measure space satisfying ${\mathscr Ric}_\phi(g) \ge - (n-1)k g$ in $\mathscr{B}_{R}$ for some  
$k \ge 0$. Let $w$ be a positive and bounded solution to $({\mathscr P})$ with $0<w \le D$.
 Then for all $(x,t)$ in $Q_{R/2,T}$ with $t>t_0-T$ we have:
\begin{align} \label{eq13}
\frac{|\nabla w|}{w} (x,t) 
\le &~C  \left(1 - \log \frac{w}{D} (x,t) \right) \times \\
& \left \{ \begin {array}{ll}  \sqrt{k}+\dfrac{1}{\sqrt{t-t_0+T}} + \dfrac{1}{R}\\
+\sqrt{\dfrac{[\gamma_{\Delta_\phi}]_+}{R}}
+ \sup_{Q_{R, T}} \left[\dfrac{|\mathscr G_x(t,x,w)|}{ w[1- \log(w/D)]^2}\right]^\frac{1}{3} \\
+  \sup_{Q_{R, T}} \left[ \dfrac{w[1-\log(w/D)] \mathscr G_w(t,x,w) 
+ \log(w/D) \mathscr G (t,x,w)}{w [1-\log(w/D)]^2} \right]_+^\frac{1}{2} \nonumber 
\end{array} 
 \right\}. 
\end{align}
Here $C>0$ is a constant depending only on the dimension $n$, $[\gamma_{\Delta_\phi}]_+=\max(\gamma_{\Delta_\phi}, 0)$ where  the functional quantity 
$\gamma_{\Delta_\phi}$ is defined by 
\begin{align} \label{sigma def}
\gamma_{\Delta_\phi} = \max_{\partial {\mathscr B}_1} \Delta_\phi r (x), 
\end{align}
with $\partial {\mathscr B}_1= \{x: d(x,x_0)=1\}$.
\end{theorem}

The local estimate above has a global counterpart subject to the prescribed bounds in the theorem 
being global. The proof follows by passing to the limit $R \to \infty$ in \eqref{eq13} and making note 
of the constants being independent of $R$ whilst $1/R$, $[\gamma_{\Delta_\phi}]_+/R \to 0$ in the 
limit. The precise formulation of the estimate is given below.

\begin{theorem} \label{thm1-global}
Let $(M, g, d\mu)$ with $d\mu=e^{-\phi} dv_g$ be a complete smooth metric measure space satisfying 
${\mathscr Ric}_\phi (g) \ge -(n-1)k g$ in $M$ for some $k \ge 0$. Let $w$ be a positive and 
bounded solution to $({\mathscr P})$ with $0<w \le D$. Then for all $x \in M$ and $t_0-T<t \le t_0$ we have:
\begin{align}\label{eq13-global}
\frac{|\nabla w|}{w}(x,t)
&\le C \left(1 - \log \frac{w}{D}(x,t) \right) \times \\
&\left \{ \begin {array}{ll}  \sqrt{k} +\dfrac{1}{\sqrt{t-t_0+T}}\\
+  \sup_{M \times [t_0-T, t_0]}\left[ \dfrac{|\mathscr G_x(t,x,w)|}{ w[1- \log(w/D)]^2}\right]^\frac{1}{3}\\
+  \sup_{M \times [t_0-T, t_0]}  \left[ \dfrac{w[1-\log(w/D)] \mathscr G_w(t,x,w) 
+ \log(w/D) \mathscr G (t,x,w)}{w [1-\log(w/D)]^2} \right]_+^\frac{1}{2}
\end{array} \right\}.\nonumber
\end{align}
Here $C>0$ is a constant depending only on the dimension $n$.
\end{theorem}

One of the useful consequences of the estimates above is the following elliptic Harnack inequality for bounded positive 
solutions to $({\mathscr P})$. Other implication will be presented subsequnetly.

\begin{theorem} \label{cor Harnack}
Under the assumptions and local bounds in Theorem $\ref{thm1}$, for every $x_1$, $x_2$ in ${\mathscr B}_{R/2}$ 
and $t_0-T< t \le t_0$ with $d=d(x_1, x_2)$, we have the Harnack inequality:
\begin{equation} \label{eq2.5-Har}
\frac{w(x_1, t)}{eD} \le \left[ \frac{w(x_2, t)}{eD} \right]^\alpha,  
\end{equation}
where the exponent $\alpha=\alpha(R) \in (0,1)$ is given explicitly by 
\begin{equation} \label{Eq-alpha-R}
{\rm exp}
\left[
 - dC \left \{ \begin {array}{ll}
\sqrt{k}+\dfrac{1}{\sqrt{t-t_0+T}} + \dfrac{1}{R} + \\ 
\sqrt{\dfrac{[\gamma_{\Delta_\phi}]_+}{R}} 
+ \sup_{Q_{R, T}} \left[\dfrac{|\mathscr G_x(t,x,w)|}{ w[1- \log(w/D)]^2}\right]^\frac{1}{3} + \\
\sup_{Q_{R, T}} \left[ \dfrac{w[1-\log(w/D)] \mathscr G_w(t,x,w) 
+ \log(w/D) \mathscr G (t,x,w)}{w [1-\log(w/D)]^2} \right]_+^\frac{1}{2} 
\end{array}
\right\}
\right].
\end{equation}
If the prescribed bounds are global as in Theorem $\ref{thm1-global}$, then for every $x_1$, $x_2$ in $M$ and $t_0-T<t \le t_0$, 
the same inequality \eqref{eq2.5-Har} holds now with the exponent $\alpha=\alpha(M)$ given by 
\begin{equation}
{\rm exp}
\left[ 
-dC  
\left \{ \begin {array}{ll}
\sqrt{k}+\dfrac{1}{\sqrt{t-t_0+T}} + 
\sup_{M \times [t_0-T, t_0]} \left[\dfrac{|\mathscr G_x(t,x,w)|}{ w[1- \log(w/D)]^2}\right]^\frac{1}{3} + \\
\sup_{M \times [t_0-T, t_0]} \left[ \dfrac{w[1-\log(w/D)] \mathscr G_w(t,x,w) 
+ \log(w/D) \mathscr G (t,x,w)}{w [1-\log(w/D)]^2} \right]_+^\frac{1}{2} 
\end{array} \right\}
\right].
\end{equation}
\end{theorem}

The second type of gradient estimates to be considered in this paper are of Hamilton type. 
Here again $x_0 \in M$, $t_0 \in {\mathbb R}$, $R>2$, $T>0$ are fixed. Unlike the Souplet-Zhang estimate in 
Theorem \ref{thm1} here there is a choice of a free parameters $\alpha>1+\beta$, $\beta \ge 0$ 
in the estimate.

\begin{theorem} \label{thm18}
Let $(M, g, d\mu)$ with $d\mu=e^{-\phi} dv_g$ be a complete smooth metric measure space satisfying 
${\mathscr Ric}_\phi (g) \ge - (n-1)k g$ in $\mathscr{B}_{R}$ for some $k \ge 0$. Let $w$ be a positive 
solution to $({\mathscr P})$. 
Then for every $\beta \ge 0$ and $\alpha>1+\beta>0$ 
and for all $(x,t)$ in $Q_{R/2,T}$ with $t>t_0-T$ we have:
\begin{align} \label{eq1.12}
\frac{|\nabla w|}{w^{1-(\beta+2)/(2\alpha)}}(x,t)
\le&~ C \Big( \sup_{Q_{R, T}}w (x,t)\Big)^{(\beta+2)/(2\alpha)} \times \\
&\left \{ \begin {array}{ll} 
\sqrt k +\dfrac{1}{\sqrt {t-t_0+T}} + \dfrac{1}{R}\\
+ \sqrt{\dfrac{[\gamma_{\Delta_\phi}]_+}{R}}
+\sup_{Q_{R, T}} \left[\dfrac{ |\mathscr G_x (t,x,w)|}{w}\right]^\frac{1}{3} \\
+\sup_{Q_{R, T}} \left[\dfrac{2w\mathscr G_w(t,x,w) -
[2-(\beta+2)/\alpha] \mathscr G(t,x,w)}{w} \right]_+^\frac{1}{2}
\end {array}
\right\}, \nonumber
\end{align}
where $C=C(\alpha, n)>0$ and 
$\gamma_{\Delta_\phi}$ is as in \eqref{sigma def}.
\end{theorem}

\begin{remark}{\em
It is useful to point out that under the stronger curvature condition ${\mathscr Ric}^m_\phi (g) \ge -(m-1)kg$ 
one can remove in both the local estimates \eqref{eq13} and \eqref{eq1.12} the functional term 
\begin{equation} \label{eq-3.8}
([\gamma_{\Delta_\phi}]_+/R)^{1/2},
\end{equation} 
that as will be seen in the proofs, arises from an application of the Wei-Wiley weighted Laplace 
comparison theorem under the lower bound on ${\mathscr Ric}_\phi (g)$. (Under the stronger 
curvature condition, the constant $C>0$ in \eqref{eq13} and \eqref{eq1.12} will depend on $m$.) 
Also note that for $m=n$ where ${\mathscr Ric}_\phi^n(g) = {\mathscr Ric}(g)$ and $\Delta_\phi=\Delta$ (here 
$\phi$ must be constant) our estimates and results on $({\mathscr P})$ are new even 
for the classical Laplace-Beltrami operator. (See also \cite{Taheri-GE-1, Taheri-GE-2, TVahGrad}.)
}
\end{remark}

Again, the local estimate above has a global counterpart, when the asserted bounds in the theorem 
are global. The proof follows by passing to the limit $R \to \infty$ in \eqref{eq1.12}. 
This is the content of the following theorem.

\begin{theorem} \label{thm18-global}
Let $(M, g, d\mu)$ with $d\mu=e^{-\phi} dv_g$ be a complete smooth metric measure space satisfying 
${\mathscr Ric}_\phi (g) \ge - (n-1)k g$
 in $M$ for some $k \ge 0$. Let $w$ be a positive solution to $({\mathscr P})$.  
Then for every  $\beta \ge 0$ and $1+\beta<\alpha$ and all $x \in M$ and $t_0-T<t \le t_0$ we have:
\begin{align} \label{eq1.12-global}
\frac{|\nabla w|}{w^{1-(\beta+2)/(2\alpha)}}(x,t)
& \le C \Big( \sup_{M \times [t_0-T, t_0]} w(x,t) \Big)^{(\beta+2)/(2\alpha)} \times \\
&\left \{ \begin {array}{ll} 
\sqrt k +\dfrac{1}{\sqrt {t-t_0+T}} \\ 
+ \sup_{M \times [t_0-T, t_0]}  \left[ \dfrac{|\mathscr G_x (t,x,w)|}{w} \right]^\frac{1}{3} \\
+\sup_{M \times [t_0-T, t_0]}  \left[\dfrac{2w\mathscr G_w(t,x,w) -
[2-(\beta+2)/\alpha] \mathscr G(t,x,w)}{w} \right]_+^\frac{1}{2}
\end {array} \right\}, \nonumber
\end{align}
where $C=C(\alpha, n)>0$.
\end{theorem}

\subsection{Applications to Liouville-type results $({\bf I})$: ${\mathscr Ric}_\phi(g) \ge 0$}

Let us now move on to discussing some applications of the gradient estimates above. Here we focus mainly 
on applications to Liouville and global constancy type results for both the parabolic and elliptic (stationary states) 
equations associated with $({\mathscr P})$. Furthermore we assume ${\mathscr G}={\mathscr G}(w)$.

Towards this end let us first suppose that $w$ is a smooth bounded positive solution to the elliptic equation 
\begin{align} \label{elliptic PDE equation}
\Delta_\phi w + {\mathscr G}(w) = \Delta w - \langle \nabla \phi, \nabla w \rangle + {\mathscr G}(w) =0.  
\end{align}
Our main goal is to present two independent elliptic Liouville-type results for positive solutions to \eqref{elliptic PDE equation}. 
Note that since a constant solution to \eqref{elliptic PDE equation} must be a zero of ${\mathscr G}$, 
if ${\mathscr G}$ has no zeros $w>0$ then a Liouville theorem can be seen as a 
non-existence result.

\begin{theorem} \label{Liouville log thm}
Let $(M, g, d\mu)$ with $d\mu= e^{-\phi} dv_g$ be a complete smooth metric measure space 
satisfying ${\mathscr Ric}_\phi(g) \ge 0$ in $M$. Let $w$ be a positive and bounded solution 
to \eqref{elliptic PDE equation} with $0<w \le D$ and assume 
$[1-\log(w/D)] w \mathscr G'(w) + \log(w/D) \mathscr G (w) \le 0$. 
Then $w$ is a constant. In particular ${\mathscr G}(w) = 0$.
\end{theorem}

The second Liouville-type result of interest as mentioned above is formulated below. The proof of both these 
theorems will be presented in Section \ref{sec6}.

\begin{theorem} \label{Liouville two thm}
Let $(M, g, d\mu)$ with $d\mu= e^{-\phi} dv_g$ be a complete smooth metric measure space 
satisfying ${\mathscr Ric}_\phi(g) \ge 0$ in $M$. Let $w$ be a positive and bounded solution to 
\eqref{elliptic PDE equation} and assume for some $\beta \ge 0$ 
and  $1+\beta < \alpha$ we have $[1-(\beta/2+1)/\alpha] {\mathscr G}(w) -w {\mathscr G}'(w) \ge 0$. 
Then $w$ is a constant. In particular ${\mathscr G}(w) = 0$.
\end{theorem}

Let us now discuss implications and some examples of the above theorems to specific nonlinearities 
of interest arising in various contexts within geometry and mathematical physics. 
For further discussion on this see the introduction to the paper.

Here we focus on nonlinearities in the ``{\it split}" form ${\mathscr G}(w) = X(w) + w^r Y(\log w)$ where 
$X=X(w) \in {\mathscr C}^1(0, \infty)$, $r$ is a real exponent and $Y=Y(s)= Y(\log w) \in {\mathscr C}^1(-\infty, \infty)$. 
The motivation behind considering such types on nonlinearities was discussed earlier in Section \ref{sec1} 
({\it see} also Theorem \ref{LiouvilleThmEx-2.9} below). The next result is an application of Theorem \ref{Liouville log thm}.

\begin{theorem} \label{Liouville one thm}
Let $(M, g, d\mu)$ with $d\mu= e^{-\phi} dv_g$ be a complete smooth metric measure space 
satisfying ${\mathscr Ric}_\phi(g) \ge 0$ in $M$. Let $w$ be a positive bounded solution to the nonlinear elliptic equation 
\begin{equation}
\Delta_\phi w + X(w) + w^r Y(\log w) =0. 
\end{equation}
Suppose $X'(w) \le 0$ and $X(w)-wX'(w) \ge 0$. 
Moreover, when $Y \not\equiv 0$, depending on $w$ satisfying $w \ge 1$ or $0<w \le 1$ everywhere 
on $M$, assume that:
\begin{align} \label{Phi ass first Liouville}
&({\mathsf Y1}):=
\begin{cases}
\mbox{$r \ge 1$, $w \ge 1$,} \\
\mbox{$Y(s) \le 0$ and}\\
\mbox{$Y'(s) \le 0$ for all $s \ge 0$},
\end{cases}\\
&or, \nonumber \\ 
&({\mathsf Y2}):= 
\begin{cases}
\mbox{$r \le \min(\gamma, 1)$, $0<w \le 1$,} \\
\mbox{$Y(s) \ge 0$, $Y'(s) \le 0$ and}\\
\mbox{$s Y'(s) \ge \gamma Y(s)$ for all $s \le 0$ and some $\gamma \ge 0$}.
\end{cases}
\end{align}
Then $w$ is a constant. In particular $X(w) + w^r Y(\log w)=0$. 
\end{theorem}

Regarding the assumptions in the theorem note that the examples $X(w) = {\mathsf A} w^p$ (with ${\mathsf A} \ge 0$ and $p \le 0$) 
and $X(w)= {\mathsf B} w^q$ (with ${\mathsf B} \le 0$ and $q \ge 1$) satisfy the assumptions on $X$. Moreover, the example
$Y(s)=|s|^\gamma$ ($\gamma>1$) satisfies ${\mathsf Y2}$, $Y(s)=-|s|^\gamma$ satisfies ${\mathsf Y1}$ and the example $Y(s)=-s^\gamma$
($\gamma \ge 1$ odd integer) satisfies both ${\mathsf Y1}$ and ${\mathsf Y2}$. (Compare these to Theorem 4.3 in \cite{Taheri-GE-2} 
and the choice $X=X(w)$ in \eqref{polies PR} below.)

As another example consider the case where the nonlinearity is a superposition 
of power-like terms in the form 
\begin{equation} \label{polies PR}
X(w) = \sum_{j=1}^N \mathsf{A}_j w^{p_j} + \sum_{j=1}^N \mathsf{B}_j w^{q_j}. 
\end{equation}
Here $\mathsf{A}_j, \mathsf{B}_j$ are constant coefficients and $p_j, q_j$ for $1 \le j \le N$ 
are real exponents. Then we can formulate the following result which improves and extends 
earlier results on Yamabe-type problems 
({\it cf.} \cite{Dung,GidSp,Taheri-GE-1,Taheri-GE-2,[Wu18]}). 
This is an application of Theorem \ref{Liouville two thm}.

\begin{theorem} \label{Liouville.2-8}
Let $(M, g, d\mu)$ with $d\mu=e^{-\phi}dv_g$ be a complete smooth metric measure space satisfying 
${\mathscr Ric}_\phi(g) \ge 0$ in $M$. Let $w$ be a positive bounded solution to the nonlinear elliptic equation   
\begin{equation} \label{elliptic PDE 2}
\Delta_\phi w + \sum_{j=1}^N \mathsf{A}_j w^{p_j} + \sum_{j=1}^N \mathsf{B}_j w^{q_j} = 0.
\end{equation}
Assume $\mathsf{A}_j \ge 0$, $\mathsf{B}_j \le 0$ and $p_j \le 1-(\beta/2+1)/\alpha$, $q_j \ge 1-(\beta/2+1)/\alpha$ 
for $1 \le j \le N$ and for some $\beta \ge 0$ and $1+\beta <\alpha$. Then $w$ must be a constant and $\mathscr G(w)=0$. 
\end{theorem}

As a final application consider ${\mathscr G}(w) = \mathsf{A} w^p + \mathsf{B} w^q + w^r Y(\log w)$ 
with real exponents $p$, $q$, $r$, constant coefficients $\mathsf{A}$, $\mathsf{B}$ and 
$Y \in \mathscr{C}^1(\mathbb R)$. Then Theorem \ref{Liouville two thm} leads to the following.

\begin{theorem} \label{LiouvilleThmEx-2.9}
Let $(M, g, d\mu)$ with $d\mu=e^{-\phi}dv_g$ be a complete smooth metric measure space satisfying 
${\mathscr Ric}_\phi(g) \ge 0$ in $M$. Let $w$ be a positive bounded solution to the nonlinear elliptic equation   
\begin{equation} \label{elliptic PDE 2}
\Delta_\phi w + \mathsf{A} w^p + \mathsf{B} w^q + w^r Y(\log w) = 0.
\end{equation}
Assume that along the solution $w$ we have the inequality $Y'+[(\beta/2+1)/\alpha+r-1] Y \le 0$ 
along with $\mathsf{A} \ge 0$, $\mathsf{B} \le 0$, $p \le 1-(\beta/2+1)/\alpha$ and $q \ge 1-(\beta/2+1)/\alpha$ 
for some $\beta \ge 0$ and $1+\beta < \alpha$. Then $w$ must be a constant and 
$\mathsf{A} w^p + \mathsf{B} w^q + w^r Y(\log w) =0$. 
\end{theorem}

Note that the assumption $Y'+[(\beta/2+1)/\alpha+r-1] Y \le 0$ in Theorem \ref{LiouvilleThmEx-2.9} is implied 
by either of $({\mathsf H1})$ or $({\mathsf H2})$ below on $Y=Y(s)= Y(\log w) \in \mathscr{C}^1(-\infty, \infty)$: 
\begin{align} \label{Phi ass second Liouville}
&({\mathsf H1}):=
\begin{cases}
\mbox{$r \ge 1-(\beta/2+1)/\alpha$, $w \ge 1$,} \\
 \mbox{$Y \le 0$ and $Y' \le 0$ for all $s \ge 0$},
\end{cases}\\
&or, \nonumber \\ 
&({\mathsf H2}):=
\begin{cases}
\mbox{$r \le 1-(\beta/2+1)/\alpha$, $0<w \le 1$,} \\
\mbox{$Y \ge 0$ and $Y' \le 0$ for all $s \le 0$}.
\end{cases}
\end{align}

We end this section with an application of the estimates above to parabolic Liouville-type theorems. 
Here by an ancient solution $w=w(x,t)$ to $({\mathscr P})$ we mean a solution defined on $M$ for 
all negative times, that is, for all $(x,t)$ with $x \in M$ and $-\infty < t <0$.

\begin{theorem} \label{ancient}
Let $(M, g, d\mu)$ with $d\mu=e^{-\phi}dv_g$ be a complete smooth metric measure space satisfying 
${\mathscr Ric}_\phi(g) \ge 0$ in $M$. Let ${\mathscr G}'(w) \le 0$, 
${\mathscr G}(w)-w{\mathscr G}'(w) \ge 0$ and ${\mathscr G}(w) \ge a$ for some $a>0$ and all $w>0$. 
Then the equation 
\begin{equation} \label{ancient-equation}
\frac{\partial w}{\partial t} - \Delta_\phi w = {\mathscr G}(w), 
\end{equation}
does not have any positive ancient solutions satisfying $w(x,t) = e^{o(\sqrt{r(x)} + \sqrt{|t|})}$.
\end{theorem}

\subsection{A parabolic Li-Yau type gradient estimate for $({\mathscr P})$}

Let us turn to our final gradient estimate for positive solutions to $({\mathscr P})$. Towards this end it is 
convenient to introduce further notation that will appear in different stages of the analysis and serve as bounds 
in the estimates. In order to describe these, for given $\mathscr{G}$ as above of class ${\mathscr C}^2$ and 
constant $\alpha$ we set,  
\begin{align}
&\mathsf{A}_{\mathscr G}^\alpha (t,x,w) = [- \alpha w{\mathscr G}_{ww}(t,x,w) + {\mathscr G}_w(t,x,w) - {w^{-1} \mathscr G}(t,x,w)]_+, \label{Eq-1.5}\\
&\mathsf{B}_{\mathscr G}^\alpha (t,x,w) = |\alpha {\mathscr G}_{xw}(t,x,w) - w^{-1} {\mathscr G}_x(t,x,w)|, \label{Eq-1.6} \\
&\mathsf{C}_{\mathscr G} (t,x,w) = [{\mathscr G}_w(t,x,w) - w^{-1}{\mathscr G}(t,x,w)]_{+},  \label{Eq-1.7} \\
&\mathsf{D}_{\mathscr G} (t,x,w) = [-w^{-1} \Delta_\phi {\mathscr G}^x(t,x,w)]_+. \label{Eq-1.8}
\end{align}
Here, as before subscripts denote partial derivatives. Moreover $\mathscr G^x: x \mapsto \mathscr G(t,x, w)$ denotes 
the function obtained by freezing the variables $t$, $w$ and viewing $\mathscr G$ as a function of $x$ only. (Thus in particular 
we speak of $\nabla {\mathscr G}^x$ and $\Delta_\phi {\mathscr G}^x$.) Having the above notation in place we now define the 
four pairs of $\gamma$-quantities associated with a given $w=w(x,t)$ ($x \in M$, $0 \le t \le T$) by writing for fixed $x_0 \in M$, 
$R>0$ and $T>0$ :
\begin{align}
&\gamma^{{\mathscr G}, \alpha}_\mathsf{A} (R) = \sup_{H_{R, T}} \mathsf{A}_{\mathscr G}^\alpha (t,x,w), 
\qquad \, \gamma^{{\mathscr G}, \alpha}_\mathsf{A} = \sup_{M \times [0, T]} \mathsf{A}_{\mathscr G}^\alpha (t,x,w), \label{eq2.2} \\
&\gamma^{{\mathscr G}, \alpha}_\mathsf{B} (R) = \sup_{H_{R, T}} \mathsf{B}^\alpha_{\mathscr G}(t,x,w), 
\qquad \, \gamma^{{\mathscr G}, \alpha}_\mathsf{B} = \sup_{M \times [0, T]} \mathsf{B}^\alpha_{\mathscr G}(t,x,w), \label{eq2.3}  \\
&\gamma^{\mathscr G}_{\mathsf C} (R) = \sup_{H_{R, T}} \mathsf{C}_{\mathscr G}(t,x,w), 
\qquad \quad \gamma^{\mathscr G}_{\mathsf C} = \sup_{M \times [0, T]} \mathsf{C}_{\mathscr G}(t,x,w), \label{eq2.4} \\
&\gamma^{\mathscr G}_{\mathsf D} (R) = \sup_{H_{R, T}} \mathsf{D}_{\mathscr G}(t,x,w), 
\qquad \quad \gamma^{\mathscr G}_{\mathsf D} = \sup_{M \times [0, T]} \mathsf{D}_{\mathscr G}(t,x,w). \label{eq2.5}
\end{align}

For the first estimate we fix $x_0 \in M$, $R>0$, $T>0$ and formulate an upper bound for the Harnack quantity 
[{\it see} the left-hand side of \eqref{1.26}] on the cylinder ${\mathscr B}_R \times (0, T]$.

\begin{theorem} \label{thm7.1}
Let $(M, g, d\mu)$ with $d\mu=e^{-\phi} dv_g$ be a complete smooth metric measure space satisfying 
${\mathscr Ric}_\phi^m(g) \ge -(m-1) k g$ in ${\mathscr B}_{2R}$ for some $m \ge n$ and $k \ge 0$. 
Let $w$ be a smooth positive solution to $({\mathscr P})$. Then for every $\alpha>1$, 
$\varepsilon \in (0, 1)$ and all $(x,t) \in H_{R,T}$ with $t>0$ we have

\begin{align}\label{1.26}
\left[ \frac{|\nabla w|^2}{\alpha w^2} - \frac{\partial_t w}{w} 
+ \frac{{\mathscr G}}{w} \right] (x,t) \le&~ 
\frac{m \alpha}{2} 
\bigg\{ \frac{1}{R^2} \left[ \frac{m c_1 ^2 \alpha ^2}{4(\alpha -1)} 
+c_2+(m-1) c_1(1+R \sqrt{k})+2c_1^2 \right] \nonumber\\
&+\frac{1}{t} + \gamma^{\mathscr G}_{\mathsf C} (2R)\bigg\}
+ \sqrt {\frac{m \alpha}{2}}\bigg\{ \frac{m \alpha [(m-1) k 
+ \gamma^{{\mathscr G}, \alpha}_\mathsf{A} (2R)/2]^2}{2(1-\varepsilon)(\alpha-1)^2}\nonumber\\
&+ \left[\frac{3^3 m [\gamma^{{\mathscr G}, \alpha}_\mathsf{B} (2R)]^{4}}{2^5 \varepsilon \alpha (\alpha-1)^2}\right]^{1/3} 
+ \gamma^{\mathscr G}_{\mathsf D} (2R) \bigg\}^{1/2}.
\end{align}
The constants $c_1$, $c_2>0$ in \eqref{1.26}
are those appearing in the bounds \eqref{9.30} in Lemma $\ref{psi lemma}$
and the  $\gamma$-quantities are as in 
\eqref{eq2.2}-\eqref{eq2.5} $($with $2R$ replacing $R$$)$.
\end{theorem}

The local estimate above has a global in space counterpart subject to the prescribed bounds in the theorem being global in space. 
The proof follows by passing $R \to \infty$ in \eqref{1.26} and taking into account the vanishing of the term
\begin{align}
 \frac{1}{R^2} \left[ \frac{m c_1 ^2 \alpha ^2}{4(\alpha -1)} 
+c_2+(m-1) c_1(1+R \sqrt{k})+2c_1^2 \right]
\end{align}
in the limit due to the constants being independent of $R$. The precise formulation 
of this global estimate is given in the theorem below.

\begin{theorem} \label{thm28-global}
Let $(M, g, d\mu)$ with $d\mu=e^{-\phi} dv_g$ be a complete smooth metric measure space 
satisfying ${\mathscr Ric}^m_\phi (g) \ge -(m-1) k g$ in $M$ for some $m \ge n$ and $k \ge 0$. 
Let $w$ be a smooth positive solution to $({\mathscr P})$. Then for every 
$\alpha>1$, $\varepsilon \in (0, 1)$ and all $x \in M$, $0<t \le T$ we have 
\begin{align} \label{1.26-global}
\left[ \frac{|\nabla w|^2}{\alpha w^2} - \frac{\partial_t w}{w} + \frac{{\mathscr G}}{w} \right] (x,t) 
\le&~ \frac{m \alpha}{2t}
+ \sqrt {\frac{m \alpha}{2}}\bigg\{ \frac{m \alpha[(m-1) k 
+ \gamma^{{\mathscr G}, \alpha}_\mathsf{A}/2]^2}{2(1-\varepsilon)(\alpha-1)^2}\nonumber\\
&+ \left[\frac{3^3 m [\gamma^{{\mathscr G}, \alpha}_\mathsf{B}]^{4}}{2^5 \varepsilon \alpha (\alpha-1)^2}\right]^{1/3} 
+ \gamma^{\mathscr G}_{\mathsf D} \bigg\}^{1/2}
+ \frac{m \alpha}{2}\gamma^{\mathscr G}_{\mathsf C}. 
\end{align}
Here the $\gamma$-quantities on the right in \eqref{1.26-global} are as given 
by \eqref{eq2.2}-\eqref{eq2.5}.  
\end{theorem}

An immediate consequence of the local and global estimates in Theorem \ref{thm7.1} and Theorem \ref{thm28-global} 
is the following parabolic Harnack inequality on the solutions.

\begin{theorem} \label{thm38}
Under the assumptions of Theorem $\ref{thm7.1}$ let $w$ be a positive solution to $({\mathscr P})$. 
Then for every $(x_1, t_1)$, $(x_2, t_2)$ in $H_{R, T}$ with $t_2>t_1>0$ and $\alpha >1$,  
\begin{align} \label{HarIneq}
w(x_2,t_2) \ge w(x_1,t_1) {\rm exp} [(t_2-t_1) {\mathsf H} - \alpha L(x_1,x_2, t_2-t_1)]  \left(\frac{t_2}{t_1}\right)^{-m \alpha/2}.  
\end{align}
Here ${\mathsf H}$ is a constant depending only on the bounds given in Theorem $\ref{thm7.1}$ 
$[$see \eqref{H-in-Harnack-Eq}-\eqref{gamma-in-Harnack-Eq}$]$ and $L$ is given by 
\begin{align}
L(x_1,x_2, t_2-t_1) = \inf_{\zeta \in \Gamma}  \left[ \frac{1}{4(t_2-t_1)} \int_{0}^{1} |\dot \zeta(t)|^2\,dt \right], 
\end{align} 
where $\Gamma  = \Gamma(x_1, x_2) = \{ \zeta \in \mathscr{C}^1( [0,1], M) : \zeta([0,1]) \subset {\mathscr B}_R, 
\zeta(0) = x_1, \zeta(1) = x_2\}$. If the bounds are global as in Theorem $\ref{thm28-global}$ then the 
Harnack inequality \eqref{HarIneq} is global too. 
\end{theorem}

\subsection{Applications to Liouville-type results $({\bf II})$: ${\mathscr Ric}^m_\phi(g) \ge 0$}

We complement the Liouville-type results given earlier by invoking the Li-Yau gradient estimates 
above which were obtained under a different curvature condition. We begin by the following theorem 
on the stationary version of $({\mathscr P})$ that is of independent interest.

\begin{theorem} \label{thm48}
Let $(M, g, d\mu)$ with $d\mu=e^{-\phi} dv_g$  be a complete smooth metric 
measure space satisfying ${\mathscr Ric}_\phi^m(g) \ge -(m-1)kg$ in $M$. 
Let $w$ be a smooth positive solution to \eqref{elliptic PDE equation}. 
Then, for every $\alpha>1$ and $\varepsilon \in (0, 1)$ the following global 
estimate holds on $M$:  
\begin{align} \label{eqL2.26}
\frac{|\nabla w|^2}{\alpha w^2} + \frac{\mathscr G{}(w)}{w}  
\le&~\frac{m \alpha}{2} 
\left[ 
\frac{(m-1) k + \gamma^{{\mathscr G}, \alpha}_\mathsf{A}/2}{(\alpha-1) \sqrt{1-\varepsilon}} 
+ \gamma^{\mathscr G}_\mathsf{C} \right].
\end{align}
\end{theorem}

\begin{theorem} \label{coroLiouville}
Let $(M, g, d\mu)$ with $d\mu=e^{-\phi} dv_g$  be a complete smooth metric measure space 
satisfying ${\mathscr Ric}_\phi^m(g) \ge 0$ in $M$. Let $w$ be a smooth positive solution to 
\eqref{elliptic PDE equation} and assume that along the solution $w$ we have 
${\mathscr G}(w) \ge 0$, ${\mathscr G}(w) - w {\mathscr G}_w(w) \ge 0$ and  
\begin{equation}
{\mathscr G}(w)-w{\mathscr G}_w(w)+\alpha w^2 {\mathscr G}_{ww}(w) \ge 0, 
\end{equation}
for some $\alpha>1$, everywhere on $M$. Then $w$ must be a constant and ${\mathscr G}(w)=0$.
\end{theorem}

Theorem \ref{coroLiouville} leads to the following conclusion on Yamabe equations. As indicated 
earlier, since a constant solution to \eqref{elliptic PDE equation} must be a zero of ${\mathscr G}$, when 
${\mathscr G}$ has no zeros $w>0$, a Liouville theorem can be seen as a non-existence result 
for positive solutions. In the following theorem this happens precisely when ${\mathsf A}_j>0$ for 
at least one $1\le j \le d$. Compare also with Theorem \ref{Liouville.2-8} and the choices of 
curvature conditions there and here.

\begin{theorem} \label{LiouvilleThmEx}
Let $(M, g, d\mu)$ with $d\mu=e^{-\phi} dv_g$ be a complete smooth metric measure space  
satisfying ${\mathscr Ric}^m_\phi(g) \ge 0$ in $M$. Let $w$ be a positive smooth solution to the equation   
\begin{equation} \label{elliptic PDE 2}
\Delta_\phi w + \sum_{j=1}^N \mathsf{A}_j w^{p_j} = 0.
\end{equation}
If $\mathsf{A}_j \ge 0$ and $p_j \le 1$ for all $1 \le j \le N$ then $w$ must be a constant and $\sum \mathsf{A}_j w^{p_j}=0$.
\end{theorem}

As another application consider a nonlinearity in the form of a superposition of 
a logarithmic and a power-like nonlinearity with constant coefficients $\mathsf{A}, \mathsf{B}$, exponent $s$ 
and $Y \in \mathscr{C}^2(\mathbb R)$ in the form 
\begin{equation}
{\mathscr G}(w) = \mathsf{A} w Y(\log w) + \mathsf{B} w^s. 
\end{equation}
The following theorem now directly results from Theorem \ref{coroLiouville}.

\begin{theorem} \label{LiouvilleThmEx-log}
Let $(M, g, d\mu)$ with $d\mu=e^{-\phi} dv_g$ be a complete smooth metric measure space satisfying 
${\mathscr Ric}^m_\phi(g) \ge 0$. Let $w$ be a positive smooth solution to the equation   
\begin{equation} \label{elliptic PDE 3}
\Delta_\phi w + \mathsf{A} w Y(\log w) + \mathsf{B} w^s = 0.
\end{equation}
Assume that along the solution $w$ we have $Y \ge 0$, $Y' \le 0$ and 
$\alpha Y''+(\alpha-1)Y'\ge0$ for some $\alpha >1$. 
Furthermore assume $\mathsf{A}, \mathsf{B} \ge 0$ and $s \le 1$. Then $w$ must be a constant 
and $\mathsf{A} w Y(\log w) + \mathsf{B} w^s=0$. 
\end{theorem}

\section{Proof of the local Souplet-Zhang estimate in Theorem \ref{thm1}} \label{sec3}

We now attend to the proof of the results formulated in Section \ref{sec2}. Throughout this 
section we confine to the proof of Theorem \ref{thm1} and all of its necessary ingredients. 
The first two subsections develop the necessary apparatus, including an evolution inequality on 
a quantity (below called $H$) built out of the solution $w$, and the last section completes the 
proof by putting these together, using a localisation argument and finally maximum principle.

\subsection{Evolution inequality for $H=|\nabla [\log (1-h)]|^2$ with $h = \log (w/D)$}

From the positive solution $w$ satisfying the bounds $0<w \le D$ we first define a non-positive 
function by putting $h=\log (w/D)$ and then $H= |\nabla h|^2/(1-h)^2$. The evolution of $H$ 
under the equation $({\mathscr P})$ is the subject of the next lemma.

\begin{lemma} \label{lem21}
Let $w$ be a positive and bounded solution to $({\mathscr P})$ with $0<w \le D$. 
Put $h = \log (w/D)$ and let $H = |\nabla h|^2/(1-h)^2$. Then $H$ satisfies the evolution equation  
\begin{align} \label{eq21}
[\Delta_\phi - \partial_t]H =&~ \frac{ 2 {\mathscr Ric}_\phi(g) (\nabla h, \nabla h)}{(1-h)^2}  
+ \frac{2h \langle \nabla h, \nabla H \rangle}{1-h} 
+ 2 (1-h) H^2 \nonumber \\
&+ 2 \left| \frac{\nabla^2 h}{1-h} + \frac{\nabla h \otimes \nabla h}{(1-h)^2} \right|^2 
- \frac{2 \langle \nabla h, \mathscr G_x(t,x,De^h) \rangle}{D e^h (1-h)^2} \nonumber \\
&- 2 H \left[ \mathscr G_w (t,x,De^h) + \frac{ h\mathscr G (t,x,De^h)}{De^h(1-h)} \right].
\end{align}
In particular, if ${\mathscr Ric}_\phi(g) \ge -(n-1) kg$ for some $k \ge 0$, then 
\begin{align} \label{eq21-Ric}
[\Delta_\phi - \partial_t]H \ge& -2(n-1)k H 
+ \frac{2h \langle \nabla h, \nabla H \rangle}{1-h} \nonumber \\
&+ 2 (1-h) H^2 
- \frac{2 \langle \nabla h, \mathscr G_x(t,x,De^h) \rangle}{D e^h (1-h)^2} \nonumber \\
&- 2 H \left[ \mathscr G_w (t,x,De^h) + \frac{ h\mathscr G (t,x,De^h)}{De^h(1-h)} \right].
\end{align}

\end{lemma}

\begin{proof}
As $w$ is a solution to $({\mathscr P})$ and $0<w \le D$ it is easily seen that the non-positive function 
$h = \log (w/D)$ satisfies the equation 
\begin{equation} \label{h evolution eq}
\partial_t h = \Delta_\phi h + |\nabla h|^2 + D^{-1} e^{-h} \mathscr G(t,x,D e^h).
\end{equation} 
As a result, by noting of $\nabla \mathscr G(t,x,De^h) = \mathscr G_x(t,x,De^h)+D e^h\nabla h \mathscr G_w (t,x,De^h) $ it follows that 
\begin{align} \label{EE-eq2.3}
\partial_t |\nabla h|^2=&~2\langle \nabla h, \nabla \partial_t h \rangle 
=2 \langle \nabla h, \nabla [\Delta_\phi h + |\nabla h|^2 + D^{-1} e^{-h} \mathscr G(t,x,De^h)] \rangle \nonumber\\
=&~ 2 \langle \nabla h, \nabla \Delta_\phi h \rangle + 2\langle \nabla h, \nabla |\nabla h|^2 \rangle 
+ \frac{2 \langle \nabla h,\mathscr G_x(t,x,De^h) \rangle}{De^h}\nonumber\\
& + 2 |\nabla h|^2 \left[\mathscr G_w(t,x,De^h) - \frac{\mathscr G(t,x,De^h)}{De^h} \right]. 
\end{align}

Now moving on to the function $H = |\nabla h|^2/(1-h)^2 = |\nabla \log (1-h)|^2$ it is easily seen that $\partial_t H 
= [\partial_t |\nabla h|^2]/(1-h)^2 + [2|\nabla h|^2 \partial_t h]/(1-h)^3$ and so making note of \eqref{h evolution eq} 
and \eqref{EE-eq2.3} we can write 
\begin{align}\label{eq22}
\partial_t H 
=&~\frac{2 \langle \nabla h, \nabla \Delta_\phi h\rangle}{(1-h)^2} 
+ \frac{2 \langle \nabla h, \nabla |\nabla h|^2 \rangle }{(1-h)^2} 
+ \frac{2 \langle \nabla h,\mathscr G_x(t,x,De^h) \rangle}{D e^f (1-h)^2} \nonumber \\
&+ \frac{2 |\nabla h|^2}{(1-h)^2} \left[\mathscr G_w(t,x,De^h) - \frac{\mathscr G(t,x,De^h)}{De^h} \right] 
 + \frac{2 |\nabla h|^2  \Delta_\phi h}{(1-h)^3} \nonumber\\
&+ \frac{2 |\nabla h|^4}{(1-h)^3} + \frac{2 |\nabla h|^2 \mathscr G(t,x,De^h)}{D e^h (1-h)^3}.
\end{align}
Likewise we have $\nabla H  = [\nabla |\nabla h|^2]/(1-h)^2 + [2 |\nabla h|^2 \nabla h]/(1-h)^3$ and so 
by recalling the formulation $\Delta_\phi H  = \Delta H - \langle \nabla \phi, \nabla H \rangle$ it follows that 
\begin{align}\label{eq23}
\Delta_\phi H 
=& \, \frac{\Delta_\phi |\nabla h|^2}{(1-h)^2} + \frac{4 \langle \nabla h, \nabla |\nabla h|^2 \rangle}{(1-h)^3} 
+ \frac{2 |\nabla h|^2 \Delta_\phi h}{(1-h)^3} + \frac{6 |\nabla h|^4}{(1-h)^4}. 
\end{align}
Putting (\ref{eq22})-(\ref{eq23}) together and taking into account the necessary cancellations give   
\begin{align}
[\Delta_\phi - \partial_t] H 
= &~\frac{\Delta_\phi |\nabla h|^2}{(1-h)^2} 
- \frac{2 \langle \nabla h, \nabla \Delta_\phi h \rangle}{(1-h)^2} 
- \frac{2 \langle \nabla h, \nabla |\nabla h|^2 \rangle}{(1-h)^2}\nonumber\\  
&- \frac{2 |\nabla h|^4}{(1-h)^3} 
+ \frac{4 \langle \nabla h, \nabla |\nabla h|^2 \rangle}{(1-h)^3} + \frac{6 |\nabla h|^4}{(1-h)^4}
 - \frac{2 \langle \nabla h, \mathscr G_x(t,x,De^h) \rangle}{D e^h (1-h)^2}\nonumber\\ 
&- \frac{2 |\nabla h|^2}{(1-h)^2} \left[\mathscr G_w(t,x,De^h) - \frac{\mathscr G(t,x,De^h)}{De^h} \right] 
- \frac{2 |\nabla h|^2 \mathscr G(t,x,De^h)}{D e^h (1-h)^3}.
\end{align}
Now by making use of the weighted Bochner-Weitzenb\"ock formula we obtain  
\begin{align}
[\Delta_\phi - \partial_t] H 
=&~\frac{2 |\nabla \nabla h|^2}{(1-h)^2} + \frac{2 {\mathscr Ric}_\phi (g) (\nabla h, \nabla h)}{(1-h)^2} 
- \frac{2 \langle \nabla h, \nabla |\nabla h|^2 \rangle}{(1-h)^2}  \nonumber \\
&+ \frac{4 \langle \nabla h, \nabla |\nabla h|^2 \rangle}{(1-h)^3} 
- \frac{2 |\nabla h|^4}{(1-h)^3}+ \frac{6 |\nabla h|^4}{(1-h)^4} 
 - \frac{2 \langle \nabla h,\mathscr G_x(t,x,De^h) \rangle}{D e^h (1-h)^2} \nonumber\\
&- \frac{2 |\nabla h|^2}{(1-h)^2} \left[ \mathscr G_w(t,x,De^h) - \frac{\mathscr G(t,x,De^h)}{De^h} \right] 
- \frac{2 |\nabla h|^2 \mathscr G(t,x,De^h)}{D e^h (1-h)^3}, 
\end{align}
and therefore a rearrangement of terms and basic considerations leads to 
\begin{align}
[\Delta_\phi - \partial_t] H 
= &~\frac{2 {\mathscr Ric}_\phi (g) (\nabla h, \nabla h)}{(1-h)^2}
+ 2 \left| \frac{\nabla \nabla h}{1-h} + \frac{\nabla h \otimes \nabla h}{(1-h)^2} \right|^2 + \frac{2 |\nabla h|^4}{(1-h)^3}\nonumber\\
&- \frac{2 \langle \nabla h, \nabla |\nabla h|^2 \rangle}{(1-h)^2}  - \frac{4 |\nabla h|^4}{(1-h)^3}
 + \frac{2 \langle \nabla h, \nabla |\nabla h|^2 \rangle}{(1-h)^3}+ \frac{4 |\nabla h|^4}{(1-h)^4}\nonumber\\
&- \frac{2 \langle \nabla h, \mathscr G_x(t,x,De^h) \rangle}{D e^h (1-h)^2}
- \frac{2 |\nabla h|^2 \mathscr G(t,x,De^h)}{D e^h (1-h)^3}\nonumber\\
&- \frac{2 |\nabla h|^2}{(1-h)^2} \left[\mathscr G_w(t,x,De^h) - \frac{\mathscr G(t,x,De^h)}{De^h} \right].
\end{align}
Finally by making note of the relation $(1-h)^3 \langle \nabla h, \nabla H \rangle=  
(1-h) \langle \nabla h, \nabla |\nabla h|^2 \rangle + 2 |\nabla h|^4$ the expression on the second line on the right simplifies to 
$2 h \langle \nabla h, \nabla H \rangle /(1-h)$ and so a reference to $H=|\nabla h|^2/(1-h)^2$ leads at once to 
the desired conclusion. 

The inequality \eqref{eq21-Ric} now follows from the above by using the Ricci curvature lower bound 
${\mathscr Ric}_\phi(g) \ge -(n-1) kg$ and the non-negativity of the quadratic term on the first line on 
the right. The proof is thus complete. 
\end{proof}

\subsection{Construction of space-time cut-offs and cylindrical localisation}

In order to prove the estimate in Theorem \ref{thm1} we first need to establish a local version of the evolution 
inequality in Lemma \ref{lem21} suitable for the application of maximum principle. This will be achieved through the use of suitable space-time cut-off functions. 
Towards this end, let us fix $x_0 \in M$, $t_0 \in {\mathbb R}$, $R, T>0$ and then $\tau \in (t_0-T, t_0]$. 
The following standard lemma grants the existence of a smooth function $\bar{\eta}=\bar \eta(r, t)$ of two real 
variables $r \ge 0$ and $t_0-T \le t \le t_0$ respectively satisfying a set of useful bounds and properties for carrying 
out this {\it cylindrical} localisation procedure (see, e.g., \cite{[Ba1], Bri, SZ}).
\begin{lemma} \label{phi lemma} Fix $t_0 \in {\mathbb R}$ and let $R, T>0$. Given $\tau \in (t_0-T, t_0]$ there exists a smooth function 
$\bar{\eta}:[0,\infty) \times [t_0-T, t_0] \to \mathbb{R}$ such that the following properties hold:
\begin{enumerate}[label=$(\roman*)$]
\item ${\rm supp} \, \bar{\eta}(r,t) \subset [0,R] \times [t_0-T, t_0]$  and $0 \leq \bar{\eta}(r,t) \leq 1$ in $[0,R] \times [t_0-T, t_0]$,
\item $\bar{\eta}=1$ in $[0,R/2] \times [\tau, t_0]$ and $\partial \bar{\eta}/\partial r =0$ in $[0,R/2] \times [t_0-T, t_0]$, respectively,
\item there exists $c>0$ such that 
\begin{equation}
\left| \frac{\partial \bar{\eta}}{\partial t} \right| \le \frac{c \bar{\eta}^{1/2}}{\tau-t_0+T}, 
\end{equation} 
in $[0,\infty)\times[t_0-T,t_0]$ and $\bar{\eta}(r,t_0-T)=0$ for all $r \in [0,\infty)$.
\item $-c_a \bar{\eta}^a/R \le \partial \bar{\eta}/\partial r \leq 0$ and $|\partial^2 \bar{\eta} / \partial r^2| \le c_a \bar{\eta}^a / R^2$ 
hold on $[0, \infty)\times [t_0-T, t_0]$ for every $0<a<1$ and some $c_a>0$.
\end{enumerate}
\end{lemma}

Having the above lemma at our disposal we now move on to introducing a smooth space-time 
cut-off function $\eta=\eta(x,t)$ by setting, for $R \ge 2$, $T>0$ and $t_0-T <\tau \le t_0$ (the 
reason for the choice $R \ge 2$ will be clear later), 
\begin{equation} \label{cut-off def}
\eta(x,t) = \bar{\eta}(r(x), t), \qquad (x, t) \in M \times [t_0-T, t_0].
\end{equation}
It is plain that $\eta$ is supported in the compact space-time cylinder $Q_{R,T} \subset M \times [t_0-T, t_0]$.  
Furthermore in view of \eqref{cut-off def} we have 
$\Delta_f \eta = \bar \eta_{r r} |\nabla r|^2 + \bar \eta_r \Delta_f r$ 
and $\partial_t \eta = \bar\eta_t$.

Let us now move on to the following useful product identity for the action of the weighted heat operator on the space-time localised function $\eta w$.

\begin{lemma} \label{product lemma} 
Let $\eta=\eta(x,t)$ be as above and let $H=H(x,t)$ be a space-time function of class $\mathscr{C}^2$. 
Then we have  
\begin{equation}
[\Delta_\phi - \partial_t](\eta H) = \eta [\Delta_\phi - \partial_t] H 
+ 2 [\langle\nabla \eta ,\nabla (\eta H) \rangle- |\nabla \eta|^2 H]/\eta 
+ H [\Delta_\phi - \partial_t] \eta.  
\end{equation}
\end{lemma}

\begin{proof}
Firstly 
$[\Delta_\phi - \partial_t] (\eta H) = \eta [\Delta_\phi - \partial_t] H 
+  2\langle \nabla H , \nabla \eta \rangle+ H [\Delta_\phi - \partial_t] \eta$ 
as is easily seen by direct differentiation. Next, 
$\langle\nabla \eta, \nabla (\eta H) \rangle
=\langle\nabla \eta, H\nabla\eta + \eta \nabla H \rangle
=H|\nabla\eta|^2 + \eta \langle \nabla\eta, \nabla H\rangle$ 
and therefore using the latter relation to substitute for the middle term 
$2 \langle \nabla H, \nabla \eta \rangle$ in the first identity gives the desired result. 
\end{proof}

\subsection{Finalising the proof of Theorem \ref{thm1}}

Let us now fix $\tau \in (t_0-T, t_0]$ and put $\eta(x,t) = \bar{\eta}(r(x), t)$ with $\bar \eta$ as in Lemma \ref{phi lemma}. 
We show that the desired estimate holds at all points $(x,\tau)$ with $d(x, x_0) \le R/2$. The arbitrariness of $\tau$ will then 
give the assertion for all $(x, t)$ in $Q_{R/2,T} $ with $t>t_0-T$. 
Now starting from the inequality \eqref{eq21-Ric} in Lemma $\ref{lem21}$ we have  
\begin{align}\label{eq2112}
[\Delta_\phi - \partial_t] H \ge & \, 
2(1-h)H^2 - 2(n-1)k H + \frac{2h \langle\nabla h,\nabla H\rangle}{1-h} \nonumber \\
&- \frac{2 \langle \nabla h, {\mathscr G}_x(t,x,De^h) \rangle}{D e^h (1-h)^2} 
 - 2 H \left[{\mathscr G}_w (t,x,De^h) + \frac{h \mathscr G (t,x,De^h)}{De^h(1-h)} \right].
\end{align}
Hence for the localised function $\eta H$ we have upon invoking Lemma \ref{product lemma} the inequality 
\begin{align} \label{eq25}
[\Delta_\phi - \partial_t](\eta H) 
\ge & \left\langle \frac{2h\nabla h}{1 -h} + \frac{2\nabla \eta}{\eta}, \nabla (\eta H) \right\rangle 
- \left\langle H \left[ \frac{2h \nabla h}{1 -h} + \frac{2 \nabla \eta}{\eta} \right], \nabla \eta \right\rangle \\
&+ 2(1-h) \eta H^2 + H [\Delta_\phi - \partial_t -2(n-1) k] \eta \nonumber \\
&- \frac{2\eta \langle \nabla h, {\mathscr G}_x(t,x,De^h) \rangle}{D e^h (1-h)^2} 
- 2 \eta H \left [{\mathscr G}_w (t,x,De^h) + \frac{h \mathscr G (t,x,De^h)}{De^h(1-h)} \right]. \nonumber
\end{align}

Assume that the localised function $\eta H$ is maximal at the point $(x_1, t_1)$ in the compact space-time 
cylinder $\{(x,t) : d(x,x_0) \le R, t_0-T \le t \le \tau\}$. Suppose also that $x_1$ is not in the cut locus of $M$ by 
Calabi's argument \cite{[LY86]} and that $(\eta H)(x_1, t_1) >0$ as otherwise the result is trivial with 
$w(x, \tau) \leq 0$ for all $d(x, x_0) \le R/2$. It then follows that $t_1>t_0-T$ at $(x_1,t_1)$ and subsequently by maximum principle
\begin{align}
\begin{cases}
\partial_t (\eta H) \ge 0, \\
\nabla (\eta H)=0, \\
\Delta (\eta H) \le 0.
\end{cases}
\end{align}
In particular at $(x_1,t_1)$ we have $[\Delta_\phi - \partial_t](\eta H) \le 0$ and therefore from \eqref{eq25} we obtain the inequality
\begin{align} 
2(1-h)\eta H^2  \leq 
& \left\langle  H \left[ \frac{2h \nabla h}{1 -h} + \frac{2 \nabla \eta}{\eta} \right], \nabla \eta \right\rangle 
- H \Delta_\phi \eta \nonumber\\
&+ H \partial_t \eta + 2(n-1)kH\eta + \frac{2 \eta \langle \nabla h, {\mathscr G}_x(t,x,De^h) \rangle}{D e^h (1-h)^2} \nonumber\\
&+ 2\eta H \left[{\mathscr G}_w (t,x,De^h) + \frac{h \mathscr G (t,x,De^h)}{De^h(1-h)} \right].
\end{align}
As a result dividing through by $2(1-h) \ge 0$ it follows that
\begin{align}\label{eq27}
\eta H^2 \le 
& \left \langle H \left[ \frac{h \nabla h}{1 -h} + \frac{\nabla \eta}{\eta} \right], \frac{\nabla \eta}{1-h} \right \rangle
+ \frac{H[-\Delta_\phi+ \partial_t+ 2(n-1)k]\eta}{2(1-h)} \nonumber \\
&+ \frac{\eta \langle \nabla h, {\mathscr G}_x(t,x,De^h) \rangle}{D e^h (1-h)^3} 
+ \eta H \left[\frac{{\mathscr G}_w (t,x,De^h)}{1-h} + \frac{h \mathscr G (t,x,De^h)}{De^h(1-h)^2} \right].
\end{align}

The goal is now to utilise \eqref{eq27} to establish the required estimate at $(x, \tau)$. Towards this end we proceed 
by considering two alternatives: $d(x_1, x_0) \le 1$ and $d(x_1, x_0) \ge 1$. \\
{\bf Alternative 1 ($d(x_1, x_0) \le 1$).} Since here $\eta$ is a constant 
function in the space [for all $x$ with $d(x, x_0) \le R/2$ and $R \geq 2$ by property 
$(ii)$] the terms involving space derivatives of $\eta$ at $(x_1, t_1)$ vanish (that is, $\nabla \eta=0$, 
$\Delta_\phi \eta =0$ and $\partial_t \eta = \bar \eta_t$). Thus \eqref{eq27} reduces to  
\begin{align}
\eta H^2 \le &~\frac{H[|\partial_t \eta|+2(n-1)k\eta]}{2(1-h)} 
+ \frac{\eta|\langle \nabla h, {\mathscr G}_x(t,x,De^h) \rangle|}{D e^h (1-h)^3} \nonumber \\
& + \eta H \left[ \frac{{\mathscr G}_w (t,x,De^h)}{1-h} + \frac{h \mathscr G (t,x,De^h)}{De^h(1-h)^2} \right]_+, 
\end{align}
and so 
\begin{align}
\eta H^2 \le&~\frac{\sqrt \eta H}{2(1-h)} \frac{|\partial_t \bar \eta|}{\sqrt{\bar \eta}} + \frac{(n-1)k \eta H}{(1-h)} 
+ \eta \sqrt H \frac{|{\mathscr G}_x(t,x,De^h)|}{D e^h (1-h)^2} \nonumber \\
&+ \eta H \left[ \frac{{\mathscr G}_w (t,x,De^h)}{1-h} + \frac{h \mathscr G (t,x,De^h)}{De^h(1-h)^2} \right]_+.
\end{align}
Next by an application of Cauchy-Schwarz and Young's inequality and after rearranging terms it follows that for suitable 
$C=C(n)>0$ we have
\begin{align}
\eta H^2 \le C \bigg\{& \frac{1}{(\tau-t_0+T)^2} + k^2 
+ \left[ \frac{|{\mathscr G}_x(t,x,De^h)|}{D e^h (1-h)^2} \right]^{4/3} \nonumber\\
&+ \left[ \frac{{\mathscr G}_w (t,x,De^h)}{1-h} + \frac{h \mathscr G (t,x,De^h)}{De^h(1-h)^2} \right]^2_+ \bigg\}.
\end{align}
As $\eta \equiv 1$, when $d(x,x_0) \le R/2$ and $\tau \le t \le t_0$ (and so in particular for when $t = \tau$) by $(i)$, we have 
$\sqrt H(x,\tau) =\sqrt{\eta H} (x, \tau) \le \sqrt{\eta H}(x_1,t_1) \le \sqrt[4]{\eta H^2} (x_1,t_1)$. 
Hence recalling $H=|\nabla h|^2/(1-h)^2$ and the relation $h=\log(w/D)$, we arrive at the bound at $(x, \tau)$ 
\begin{align}
\frac{\left| \nabla h \right|}{1- h} \le C \bigg\{&\frac{1}{\sqrt {\tau-t_0+T}} + \sqrt k 
+ \sup_{Q_{R, T}} \bigg[ \frac{|{\mathscr G_x(t,x,De^h)}|}{De^h (1-h)^2} \bigg]^{1/3}\nonumber\\
&+ \sup_{Q_{R, T}} \bigg[ \frac{De^h (1-h) {\mathscr G}_w(t,x,De^h) + h \mathscr G(t,x,De^h)}{De^h(1-h)^2} \bigg]_+^{1/2} \bigg\}.   
\end{align}
This together with the arbitrariness of $\tau>t_0-T$ is now seen to give a special case of \eqref{eq13}.

\qquad \\
{\bf Alternative 2 ($d(x_1, x_0) \ge 1$).} Upon referring to the right-hand side of \eqref{eq27}, and noting the properties 
of $\bar \eta$ as listed in Lemma \ref{phi lemma} we proceed onto bounding each of the 
four terms on the right-hand side on \eqref{eq27}.

\begin{itemize}
\item Towards this end dealing with the first term first, we have 
\begin{align} \label{eq28}
\left \langle H \left[ \frac{h \nabla h}{1 -h} + \frac{\nabla \eta}{\eta} \right], \frac{\nabla \eta}{1-h} \right \rangle 
&\le H \left[ \frac{h |\nabla h|}{1 -h} + \frac{|\nabla \eta|}{\eta} \right] \frac{|\nabla \eta|}{1-h} \\
&\le H \left[ h \sqrt H + \frac{|\nabla \eta|}{\eta} \right] \frac{|\nabla \eta|}{1-h} \nonumber \\
&\le H \sqrt \eta \left [ \eta^{1/4} \frac{\sqrt H |h|}{1-h} \frac{|\nabla \eta|}{\eta^{3/4}} 
+ \frac{|\nabla \eta|^2}{\eta^{3/2}} \right] \nonumber \\
&\le \frac{\eta H^2}{5} + \frac{C}{R^4}. \nonumber 
\end{align}
\item Here we utilise the Wei-Wylie weighted Laplacian comparison theorem along with 
the lower curvature bound on ${\mathscr Ric}_\phi(g)$. 
Indeed recalling $r(x)=d(x, x_0)$, $1 \le d(x_1,x_0) \leq R$ 
and ${\mathscr Ric}_\phi (g) \ge -(n-1)k g$ with $k \ge 0$, 
upon referring to \eqref{cut-off def} we can write 
\begin{align}
-\Delta_\phi \eta &= -(\bar \eta_{r r} |\nabla r|^2 + \bar \eta_r \Delta_\phi r) \nonumber \\
&\le -(\bar \eta_{r r} + \bar \eta_r [\gamma_{\Delta_\phi} +(n-1)k (R-1)])
\end{align}
where we have used Theorem 3.1 in \cite{[WeW09]} to write 
$\Delta_\phi r \le [\gamma_{\Delta_\phi}]_+ +(n-1)k (R-1)$ 
whenever $1 \le r \le R$, $t_0-T \le t \le t_0$ [thus in particular at the space-time point $(x_1, t_1)$]. 
Note that here we have used $(ii)$ [$\bar \eta_r =0$ 
when $0 \le r \le R/2$] and $(iv)$ [$\bar \eta_r \le 0$ when $0 \le r < \infty$] 
in Lemma \ref{phi lemma}. Hence putting the above together, by an application of Young's 
inequality, we can write 
\begin{align} \label{bound Delta f alpha eq}
\frac{H (-\Delta_\phi \eta)}{1-h} &
\le \frac{\sqrt{\eta} H}{1-h}
\left[ \frac{|\bar \eta_{rr}|}{\sqrt{\bar \eta}} + \frac{|\bar \eta_r|}{\sqrt{\bar \eta}} 
\left[ [\gamma_{\Delta_\phi}]_+ +(n-1)k (R-1) \right] \right]  \nonumber\\
&\le C \frac{\sqrt \eta H}{1-h} \left[\frac{1}{R^2} + \frac{[\gamma_{\Delta_\phi}]_+}{R} + k \right]  \nonumber \\
&\le \frac{\eta H^2}{10} + C \left[ \frac{1}{R^4} 
+ \frac{[\gamma_{\Delta_\phi}]_+^2}{R^2} + k^2 \right],
\end{align}
and in a similar way using $(iii)$ in Lemma \ref{phi lemma}, 
$\partial_t \eta = \bar{\eta}_t \le c\sqrt{\eta}/(\tau-t_0+T)$, and therefore
$[2(n-1)k + \partial_t] \eta \le \sqrt \eta [2(n-1)k + c/(\tau-t_0+T)]$. An application of Young's inequality 
now gives$[H/(1-h)] [2(n-1)k + \partial_t] \eta \le \eta H^2/10 + C[1/(\tau-t_0+T)^2 + k^2]$. 
Hence returning the second term on the right in \eqref{eq27} and combining the above we have
\begin{align}
\frac{H[-\Delta_\phi+ \partial_t+ 2(n-1)k]\eta}{2(1-h)} \le~ \frac{\eta H^2}{5}
&+ C \bigg[ \frac{1}{R^4} 
+ \frac{[\gamma_{\Delta_\phi}]_+^2}{R^2} \nonumber\\
&+\frac{1}{(\tau-t_0+T)^2}+ k^2 \bigg].
\end{align}
\item In much the same way regarding the third term involving $\mathscr G (t,x,De^h)$ we have 
\begin{align}
\frac{\eta \langle \nabla h, {\mathscr G}_x(t,x,De^h) \rangle}{D e^h (1-h)^3} 
&\le \frac{\eta |\nabla h| |{\mathscr G}_x(t,x,De^h)|}{D e^h (1-h)^3} \nonumber\\
&= \frac{\eta \sqrt H |{\mathscr G}_x(t,x,De^h)|}{D e^h (1-h)^2} \nonumber \\
&\le \frac{\eta H^2}{5} + C \left[\frac{|{\mathscr G}_x(t,x,De^h)|}{De^h (1-h)^2} \right]^{4/3}.
\end{align}
\item Likewise for the subsequent terms, upon noting $-1 \le h/(1-h) \le 0$, $h \le 0$ and $0 \le \eta \le 1$, we have
\begin{align}\label{eq6.11}
\eta H  \bigg[\frac{{\mathscr G}_w(t,x,De^h)}{1-h} &+ \frac{h \mathscr G(t,x,De^h)}{De^h(1-h)^2} \bigg] \nonumber \\
=&~ \eta H \left[ \frac{De^h (1-h) {\mathscr G}_w(t,x,De^h) + h \mathscr G(t,x,De^h)}{De^h(1-h)^2} \right] \nonumber\\
\le&~ \frac{\eta H^2}{5} + C \left[ \frac{De^h (1-h) {\mathscr G}_w (t,x,De^h)
+ h \mathscr G(t,x,De^h)}{De^h (1-h)^2} \right]^2_+.
\end{align}
\end{itemize}

Now referring to \eqref{eq27} noting the inequality $1-h \ge 1$ 
and making use of the bounds obtained in \eqref{eq28}--\eqref{eq6.11}, 
it follows  the following upper bound holds for $\eta H^2$ at $(x_1, t_1)$,  
\begin{align}
\eta H^2 
\le&~ C \left \{ \begin {array}{ll} 
k^2+\dfrac{1}{(\tau-t_0+T)^2} +\dfrac{1}{R^4} \\ 
+\dfrac{[\gamma_{\Delta_\phi}]_+^2}{R^2} 
+\sup_{Q_{R, T}} \left[\dfrac{|{\mathscr G}_x(t,x,De^h)|}{De^h (1-h)^2} \right]^{4/3}\\
+\sup_{Q_{R, T}} \left[ \dfrac{De^h (1-h) {\mathscr G}_w (t,x,De^h)
+ h \mathscr G(t,x,De^h)}{De^h (1-h)^2} \right]^2_+
\end{array}
\right\}.
\end{align}

Recalling the maximality of $\eta H$ at $(x_1, t_1)$ along with $\eta \equiv 1$ when $d(x, x_0) \le R/2$ and $\tau \le t \le t_0$, 
it follows that $H^2 (x, \tau) = (\eta^2 H^2)(x, \tau) \le (\eta^2 H^2)(x_1,t_1) \leq (\eta H^2)(x_1,t_1)$ 
when $d(x, x_0) \le R/2$. hence upon noting $H = |\nabla h|^2/(1-h)^2$, the above gives
\begin{align}\label{eq214}
\frac{\left| \nabla \log w \right|}{1- \log(w/D)}
\le C \left \{ \begin {array}{ll} 
\sqrt{k} +\dfrac{1}{\sqrt{\tau-t_0+T}} + \dfrac{1}{R} \\
+\sqrt{\dfrac{[\gamma_{\Delta_\phi}]_+}{R}} 
+ \sup_{Q_{R, T}} \left[\dfrac{|{\mathscr G}_x(t,x,De^h)|}{De^h (1-h)^2} \right]^{1/3}\\
+ \sup_{Q_{R, T}} \left[\dfrac{De^h (1-h) {\mathscr G}_w (t,x,De^h)+ h \mathscr G(t,x,De^h)}{De^h (1-h)^2} \right]^{1/2}
\end{array}
\right\}.
\end{align}
Thus in either case we have shown the estimate is true at $(x, \tau)$. The desired 
conclusion now follows from the arbitrariness of $\tau \in (t_0-T, t_0]$. The proof is thus complete. \hfill $\square$

\section{Proof of the elliptic Harnack inequality in Theorem \ref{cor Harnack}} \label{sec4}

We now come to the proof of the Harnack inequality in Theorem \ref{cor Harnack}. 
To this end pick $x_1, x_2$ in $M$ and $t_0-T<t<t_0$. Let $\zeta=\zeta(s)$ with $0 \le s \le 1$ be a shortest geodesic 
curve with respect to $g$, lying completely in the closed ball ${\mathscr B}_{R/2} \subset M$ and joining the points $x_1$, 
$x_2$, specifically, $\zeta(0)=x_1$ and $\zeta(1)=x_2$. Let us also put $d=d(x_1, x_2)$.

Henceforth we shall assume the $\mathscr G$-terms in the theorem are finite 
(else $\alpha(R)=0$ and the desired inequality is trivially true). Now utilising the estimate \eqref{eq13} in Theorem \ref{thm1} 
we can write 
\begin{align} \label{Eq-5.1}
\log \frac{1-h(x_2, t)}{1-h(x_1, t)} &= \int_0^1 \frac{d}{ds} \log [1-h(\zeta(s), t)] \, ds 
= \int_0^1 - \frac{\langle \nabla h (\zeta(s), t), \zeta'(s) \rangle}{1-h(\zeta(s), t)} \, ds \\
&\le \int_0^1 \frac{|\nabla h||\zeta'|}{1-h} \, ds 
\le \sup_{Q_{R/2, T}} \left[ \frac{|\nabla h|}{1-h} \right] \int_0^1 |\zeta'|\, ds 
= d \sup_{Q_{R/2, T}} \left[ \frac{|\nabla h|}{1-h} \right] \nonumber \\
&\le dC \left\{\begin {array}{ll}  
\sqrt{k} +\dfrac{1}{\sqrt{t-t_0+T}} + \dfrac{1}{R} \nonumber \\
+\sqrt{\dfrac{[\gamma_{\Delta_\phi}]_+}{R}} 
+ \sup_{Q_{R, T}} \left[\dfrac{|{\mathscr G}_x(t,x,De^h)|}{De^h (1-h)^2} \right]^{1/3}\\
+ \sup_{Q_{R, T}} \left[\dfrac{De^h (1-h) {\mathscr G}_w (t,x,De^h)+ h \mathscr G(t,x,De^h)}{De^h (1-h)^2} \right]_+^{1/2}
\end{array}
\right\}, \nonumber
\end{align}
where in the language of \eqref{Eq-alpha-R} the last expression on the right is $-\log \alpha(R)$. 
Moreover here we have used $|\nabla h|/(1-h) = |\nabla \log w|/[1-\log(w/D)]$. 
Therefore exponentiating \eqref{Eq-5.1} results in  
\begin{align}
\frac{\log [eD/w(x_2, t)]}{\log [eD/w(x_1, t)]} 
&= \frac{1-h(x_2, t)}{1-h(x_1, t)} 
= {\rm exp} \left[ 
 \int_0^1 \frac{d}{ds} \log [1-h(\zeta(s), t)] \, ds
\right] \nonumber \\
&\le {\rm exp} [-\log \alpha(R)] = \alpha^{-1}(R), 
\end{align}
and so the claim follows at once by a further exponentiation and rearranging terms. \hfill $\square$

\section{Proof of the local Hamilton estimate in Theorem \ref{thm18}} \label{sec5}

This section is devoted to the proof of Theorem \ref{thm18} and its ingredients. The first 
subsection develops an evolution inequality on a quantity (below called $F^\alpha_\beta$) 
built out of the solution $w$, where unlike the previous case, there are now two parameters 
$\alpha>1$, $\beta \ge 0$ present, and the last section completes the proof by putting these 
together, using a localisation argument and finally concluding by an application of the maximum principle.

\subsection{Evolution inequality for  $F^\alpha_\beta = w^{(\beta+2)/\alpha-2}|\nabla w|^2/\alpha^2$}

From the positive solution $w$ we first define a positive function through $f=w^{1/\alpha}$ (with $\alpha>1$ fixed) 
and then $F^\alpha_\beta = f^\beta |\nabla f|^2$ (with $\beta \ge 0$ fixed). The evolution of $F=F^\alpha_\beta$ under 
the equation $({\mathscr P})$ is the subject of the next lemma.

\begin{lemma} \label{lemma h two}
Let $w$ be a positive solution to $({\mathscr P})$. For $\alpha>1$ and $\beta\ge 0$ fixed constants put $f=w^{1/\alpha}$ 
and $F^\alpha_\beta=f^\beta |\nabla f|^2 = w^{(\beta+2)/\alpha-2}|\nabla w|^2/\alpha^2$. 
Then $F= F^\alpha_\beta(x,t)$ satisfies the evolution equation
\begin{align} \label{eq-3.10}
[\Delta_\phi - \partial_t] F =& ~2f^\beta {\mathscr Ric}_\phi(\nabla f, \nabla f)
 + 2[1-\alpha+\beta] f^{\beta-1} \langle \nabla f, \nabla |\nabla f|^2 \rangle  \nonumber\\
 &+2f^\beta |\nabla \nabla f|^2 - [2-\beta^2 -\alpha (2-\beta)] F^2/ f^{\beta+2}\nonumber\\
&-\{2 f^{\alpha}{\mathscr G}_w(t,x, f^\alpha) -[2-(2+\beta)/\alpha]\mathscr G(t,x, f^\alpha)\} F/f^{\alpha} \nonumber\\
&- (2/\alpha) f^{\beta+1-\alpha}\langle \nabla f, {\mathscr G}_x(t,x, f^\alpha) \rangle.
\end{align}
In particular, if ${\mathscr Ric}_\phi(g) \ge -(n-1) kg$ for some $k \ge 0$, then 
\begin{align} \label{heatfW thm18}
[\Delta_\phi - \partial_t] F
\ge&- 2(n-1)k F+2(1-\alpha) \langle \nabla f, \nabla F \rangle/f \nonumber\\
&+ [2\alpha -2 -\beta^2 -\beta (2-\alpha)] F^2/f^{\beta+2} \nonumber\\
&-\{2 f^{\alpha}{\mathscr G}_w(t,x, f^\alpha)-[2-(2+\beta)/\alpha]\mathscr G(t,x, f^\alpha)\} F/ f^\alpha\nonumber\\
& - (2/\alpha) f^{\beta+1-\alpha}\langle \nabla f, {\mathscr G}_x(t,x, f^\alpha) \rangle.
\end{align}
\end{lemma}

\begin{proof}
It follows from the defining relation $f=w^{1/\alpha}$ that 
$\partial_t w = \alpha f^{\alpha-1} \partial_tf$ and 
$\Delta_\phi w = \alpha f^{\alpha-1} \Delta_\phi f + \alpha (\alpha-1) f^{\alpha-2} |\nabla f|^2$. 
Thus, from $({\mathscr P})$ it follows that 
\begin{align} \label{eq4.1}
\partial_t f = \Delta_\phi f + (\alpha -1) |\nabla f|^2/f + (1/\alpha)f^{1-\alpha} \mathscr G(t,x, f^\alpha).
\end{align}

Next, an application of the weighted Bochner-Weitzenb\"ock formula and 
the evolution of $f$ described by \eqref{eq4.1} gives, 
\begin{align}\label{eq7.4}
[\Delta_\phi - \partial_t] |\nabla f|^2  
=& \, 2 |\nabla \nabla f|^2 + 2 \langle \nabla f, \nabla \Delta_\phi f \rangle 
+ 2 {\mathscr Ric}_\phi(\nabla f, \nabla f) - \partial_t |\nabla f|^2 \nonumber \\
=& \, 2 |\nabla \nabla f|^2 + 2 \langle \nabla f, \nabla (\Delta_\phi f -\partial_t f) \rangle 
+ 2 {\mathscr Ric}_\phi(\nabla f, \nabla f) \nonumber \\
=& \, 2 |\nabla \nabla f|^2
+ 2(1-\alpha) \langle \nabla f, \nabla [|\nabla f|^2/f] \rangle \nonumber\\
&-(2 /\alpha) \langle \nabla f, \nabla [f^{1-\alpha} \mathscr G (t,x,f^\alpha)] \rangle 
+ 2 {\mathscr Ric}_\phi(\nabla f, \nabla f). 
\end{align}
In a similar way we have  
\begin{align} 
[\Delta_\phi - \partial_t] f^{\beta} 
&= \beta f^{\beta -1} [\Delta_\phi -\partial_t] f +\beta (\beta-1)f^{\beta-2} |\nabla f|^2 \nonumber \\
&= \beta f^{\beta -1}
[(1-\alpha) |\nabla f|^2/f - (1/\alpha) f^{1-\alpha} \mathscr G(t,x, f^\alpha)] 
+ \beta (\beta-1)f^{\beta-2} |\nabla f|^2 \nonumber \\
&= \beta (\beta-\alpha)f^{\beta-2} |\nabla f|^2 - (\beta/\alpha) f^{\beta-\alpha} \mathscr G(t,x, f^\alpha).
\end{align} 
As by a straightforward calculation 
\begin{align}\label{eq7.5}
[\Delta_\phi - \partial_t] (f^\beta|\nabla f|^2) 
= f^\beta [\Delta_\phi - \partial_t]  |\nabla f|^2 + 2\langle \nabla f^\beta , \nabla |\nabla f|^2 \rangle
+ |\nabla f|^2[\Delta_\phi - \partial_t] f^\beta,    
\end{align}
it then follows by putting \eqref{eq7.4}-\eqref{eq7.5} together and rearranging the equation by moving the terms 
involving the nonlinearity ${\mathscr G}$ to the end that 
\begin{align}\label{eq3.15}
[\Delta_\phi - \partial_t] F =& ~2f^\beta |\nabla \nabla f|^2 
+ 2(1-\alpha) f^{\beta -1}[\langle \nabla f, \nabla |\nabla f|^2 \rangle - |\nabla f|^4/f] \nonumber\\
&+2f^\beta {\mathscr Ric}_\phi(\nabla f, \nabla f)
+ 2 \beta f^{\beta-1}\langle \nabla f, \nabla |\nabla f|^2 \rangle\nonumber \\
&+ \beta(\beta -\alpha) f^{\beta -2} |\nabla f|^4
- (\beta/\alpha) f^{\beta-\alpha}|\nabla f|^2 \mathscr G(t,x, f^{\alpha})\nonumber\\
&-(2 /\alpha) f^\beta \langle \nabla f, \nabla [f^{1-\alpha} \mathscr G (t,x,f^{\alpha})]  \rangle.
\end{align}

Next, by virtue of $\nabla \mathscr G(t,x, f^\alpha) ={\mathscr G}_x(t,x, f^\alpha) 
+\alpha f^{\alpha -1} {\mathscr G}_w(t,x, f^\alpha) \nabla f$ we can compute 
\begin{align} \label{Pf thm18 eq3}
\langle \nabla f, \nabla [f^{1-\alpha} \mathscr G(t,x, f^\alpha)] \rangle
=&~(1-\alpha) |\nabla f|^2 \mathscr G(t,x, f^\alpha)/f^\alpha
+ f^{1-\alpha} \langle \nabla f, \nabla \mathscr G(t,x, f^\alpha) \rangle\nonumber\\
=&~ (1-\alpha) |\nabla f|^2 \mathscr G(t,x, f^\alpha)/ f^\alpha
+ \alpha|\nabla f|^2 {\mathscr G}_w(t,x, f^\alpha)\nonumber\\
& + f^{1-\alpha}\langle \nabla f, \mathscr G_x(t,x, f^\alpha) \rangle. 
\end{align}
Therefore by substituting \eqref{Pf thm18 eq3} back into the \eqref {eq3.15}, rearranging terms and recalling the relation 
$F=f^\beta|\nabla f|^2$ we can write 
\begin{align}
[\Delta_\phi - \partial_t] F 
= &~2f^\beta |\nabla \nabla f|^2 +2(1-\alpha) f^{\beta-1}\langle \nabla f, \nabla |\nabla f|^2 \rangle 
+2f^\beta {\mathscr Ric}_\phi(\nabla f, \nabla f) \nonumber\\
&+ 2 \beta f^{\beta-1}\langle \nabla f, \nabla |\nabla f|^2 \rangle
- [2-\beta^2 - \alpha(2-\beta)]f^{\beta-2} |\nabla f|^4 \nonumber\\
&-(\beta/\alpha) f^{\beta-\alpha} |\nabla f|^2 \mathscr G(t,x, f^\alpha) -2f^\beta |\nabla f|^2{\mathscr G}_w(t,x, f^\alpha)\nonumber\\
&-(2/\alpha) f^{\beta+1-\alpha}\langle \nabla f, {\mathscr G}_x(t,x, f^\alpha) \rangle
- 2[(1-\alpha)/\alpha] f^{\beta-\alpha} |\nabla f|^2 \mathscr G(t,x, f^\alpha) \nonumber\\
=&~2f^\beta |\nabla \nabla f|^2 +2(1-\alpha)f^{\beta-1} \langle \nabla f, \nabla |\nabla f|^2 \rangle 
+2f^\beta {\mathscr Ric}_\phi(\nabla f, \nabla f) \nonumber\\
&+ 2 \beta f^{\beta-1}\langle \nabla f, \nabla |\nabla f|^2 \rangle 
- [2-\beta^2 -  \alpha(2-\beta)]F^2/f^{\beta+2}  \nonumber \\
&-\{2 f^{\alpha}{\mathscr G}_w(t,x, f^\alpha) -[2-(2+\beta)/\alpha]\mathscr G(t,x, f^\alpha)\} F/f^\alpha\nonumber\\
&- (2/\alpha) f^{\beta+1-\alpha}\langle \nabla f, {\mathscr G}_x(t,x, f^\alpha) \rangle,  
\end{align}
which is the required conclusion as stated in \eqref{eq-3.10}.

To justify the inequality in the second part we make use of \eqref{eq-3.10} together with the Ricci curvature lower bound 
${\mathscr Ric}_\phi (g) \ge -(n-1)kg$ and the basic inequaltiy 
\begin{equation*}
f^\beta |\nabla \nabla f|^2 + \beta f^{\beta -1}\langle \nabla f, \nabla |\nabla f|^2 \rangle + \beta^2 f^{\beta -2} |\nabla f|^4 \ge 0,
\end{equation*}
resulting in turn from ``completing the square" ideneity
\begin{align*}
f |\nabla \nabla f|^2 + \beta \langle \nabla f, \nabla |\nabla f|^2 \rangle 
= |\sqrt{f} \nabla\nabla f + \beta[\nabla f \otimes \nabla f]/\sqrt{f}|^2 - \beta^2 |\nabla f|^4/f. 
\end{align*}
Therefore substituting in  \eqref{eq-3.10} and again recalling $F=f^\beta|\nabla f|^2$ it follows that 
\begin{align} \label{eq.3-18}
[\Delta_\phi - \partial_t] F \ge & - 2(n-1)k f^\beta |\nabla f|^2+2(1-\alpha) f^{\beta -1}\langle \nabla f, \nabla |\nabla f|^2 \rangle\nonumber\\
&- [2-\beta^2 - \alpha(2-\beta)] F^2/ f^{\beta+2} -2\beta^2 f^{\beta -2} |\nabla f|^4\nonumber\\
&-\{2 f^{\alpha}{\mathscr G}_w(t,x, f^\alpha) -[2- (2+\beta)/\alpha]\mathscr G(t,x, f^\alpha)\} F/f^\alpha \nonumber\\
&- (2/\alpha) f^{\beta+1-\alpha}\langle \nabla f, {\mathscr G}_x(t,x, f^\alpha) \rangle\nonumber\\
\ge &- 2(n-1)k F +2(1-\alpha) f^{\beta -1}\langle \nabla f, \nabla |\nabla f|^2 \rangle\nonumber\\
&- [2+\beta^2 - \alpha (2-\beta)] F^2/ f^{\beta+2} 
- (2/\alpha) f^{\beta+1-\alpha}\langle \nabla f, {\mathscr G}_x(t,x, f^\alpha) \rangle\nonumber\\
&-\{2 f^{\alpha}{\mathscr G}_w(t,x, f^\alpha) -[2-(2+\beta)/\alpha]\mathscr G(t,x, f^\alpha)\} F/f^\alpha.
\end{align}

By making note of $\langle \nabla f, \nabla F \rangle = \langle \nabla f, \nabla (f^\beta |\nabla f|^2) \rangle
 = \beta f^{\beta-1}|\nabla f|^4 + f^\beta \langle \nabla f, \nabla |\nabla f|^2 \rangle$ 
and substituting back in \eqref{eq.3-18} for $F$ we thus conclude   
\begin{align}
[\Delta_\phi - \partial_t] F
\ge&- 2(n-1)k F +2(1-\alpha) \langle \nabla f, \nabla F \rangle/f  
- 2 \beta (1-\alpha) F^2/ f^{\beta+2} \nonumber\\
&- [2+\beta^2 - \alpha(2-\beta)] F^2/ f^{\beta+2} - (2/\alpha) f^{\beta+1-\alpha}\langle \nabla f, {\mathscr G}_x(t,x, f^\alpha) \rangle\nonumber\\
&-\{2 f^{\alpha}{\mathscr G}_w(t,x, f^\alpha) -[2-(2+\beta)/\alpha]\mathscr G(t,x, f^\alpha)\}F / f^\alpha \nonumber\\
\ge&- 2(n-1)k F+2(1-\alpha) \langle \nabla f, \nabla F \rangle/f \nonumber\\
&+ [2\alpha -2 -\beta^2 - \beta(2-\alpha)] F^2/f^{\beta+2} 
- (2/\alpha) f^{\beta+1-\alpha}\langle \nabla f, {\mathscr G}_x(t,x, f^\alpha) \rangle\nonumber\\
&-\{2 f^{\alpha}{\mathscr G}_w(t,x, f^\alpha)-[2-(2+\beta)/\alpha]\mathscr G(t,x, f^\alpha)\} F / f^\alpha.
\end{align}
which is the required conclusion as stated in \eqref{heatfW thm18}. This completes the proof. 
\end{proof}

\subsection{Finalising the proof of Theorem \ref{thm18}}

Let us now turn on to completing the proof of the local estimate in Theorem \ref{thm18}. Here, for reasons that will become 
clear shortly, the ranges of the parameters $\alpha, \beta$ will be restricted to $\beta \ge 0$, $\alpha>1+\beta$ (see the last stage of the proof). 
Towards this end recalling the inequality \eqref{heatfW thm18} in Lemma \ref{lemma h two} we have 
\begin{align}
[\Delta_\phi - \partial_t] F
\ge&~2(1-\alpha) \langle \nabla f, \nabla F \rangle/f + [2\alpha -2 -\beta^2 -\beta (2-\alpha)] F^2/ f^{\beta+2} \nonumber\\
&- 2(n-1)k F - (2/\alpha) f^{\beta+1-\alpha}\langle \nabla f, {\mathscr G}_x(t,x, f^\alpha) \rangle\nonumber\\
&-\{2 f^{\alpha}{\mathscr G}_w(t,x, f^\alpha)-[2-(2+\beta)/\alpha]\mathscr G(t,x, f^\alpha)\} F/f^\alpha.
\end{align}

Next, localising by taking a space-time cut-off function $\eta$ as in \eqref{cut-off def} and following similar principles to those 
used in the proof of Theorem \ref{thm1}, we can write 
\begin{align}  \label{Pf thm18 eq5}
[\Delta_\phi - \partial_t ] (\eta F)
\ge & ~2 \langle (1-\alpha) \nabla f/f + \nabla \eta /\eta, \nabla(\eta F ) \rangle \nonumber\\
&- 2 F \langle (1-\alpha)\nabla f/f + \nabla \eta/\eta, \nabla \eta \rangle \nonumber \\
& + [2\alpha -2 -\beta^2 -\beta(2-\alpha)] \eta F^2/f^{\beta+2} \nonumber\\
&+ F [\Delta_\phi - \partial_t - 2(n-1)k]\eta \nonumber\\
&- (2/\alpha) f^{\beta+1-\alpha}\eta \langle \nabla f, {\mathscr G}_x(t,x, f^\alpha) \rangle\nonumber\\
&-\{2 f^{\alpha}{\mathscr G}_w(t,x, f^\alpha) -[2-(2+\beta)/\alpha]\mathscr G(t,x, f^\alpha)\} \eta F/f^{\alpha}.
\end{align}
For fixed $\tau \in (t_0-T,t_0]$ let $(x_1, t_1)$ be a maximum point for the localised function $\eta F$ 
in the compact set $\{(x,t) : d(x, x_0) \le R, t_0-T \le t \le \tau\}$. Then without loss of generality we 
can take $t_1>t_0-T$ and for the sake of establishing the estimate at $(x, \tau)$ in $Q_{R/2,T}$ it 
suffices to confine to the case $d(x_1, x_0) \ge 1$. Now at $(x_1, t_1)$ we have the 
inequalities
\begin{align}
\begin{cases}
\partial_t (\eta F) \ge 0, \\
\nabla (\eta F)=0, \\
\Delta (\eta F) \le 0.
\end{cases}
\end{align}
Therefore applying these to \eqref{Pf thm18 eq5} and rearranging the inequality result in 
\begin{align} \label{eq4.2}
[2\alpha -2 -\beta^2 -\beta(2-\alpha)] \eta F^2
\le &~ 2 f^{\beta+2} F \langle [ (1-\alpha) \nabla f/f + \nabla \eta/\eta], \nabla \eta \rangle \nonumber\\
&- f^{\beta+2} F [\Delta_\phi - \partial_t -2(n-1)k]\eta \nonumber\\
&+ (2/\alpha) f^{2\beta+3-\alpha} \eta \langle \nabla f, {\mathscr G}_x(t,x, f^\alpha) \rangle\nonumber\\ 
&+\{2 f^{\alpha}{\mathscr G}_w(t,x, f^\alpha) \nonumber\\
&-[2-(2+\beta)/\alpha]\mathscr G(t,x, f^\alpha)\} f^{\beta+2-\alpha}\eta F.  
\end{align}

We now proceed onto bounding from above each of the terms on the right-hand side of \eqref{eq4.2}. Again the argument proceeds by considering the 
two case $d(x_1, x_0) \le 1$ and $d(x_1, x_0) \ge 1$. However in view of certain similarities with the proof of Theorem \ref{thm1}, we shall 
remain brief, focusing on case two only and mainly on the differences.

\begin{itemize}
\item First the sum on the first line on the right-hand side of \eqref{eq4.2} is 
bounded directly in modulus by using the Cauchy-Schwarz followed by Young's inequality:  
\begin{align}\label{eq.5.15}
2 (1-\alpha) f^{\beta+1} F \langle\nabla f, \nabla \eta \rangle 
&\le 2 \left| (1-\alpha) f^{\beta+1} F \langle \nabla f, \nabla \eta \rangle \right| \nonumber \\
&\le 2(\alpha-1)\eta^{3/4} F^{3/2} f^{1+\beta/2} \frac{|\nabla \eta|}{\eta^{3/4}} \nonumber\\
&\le \frac{(\alpha-1)}{4} \eta F^2 + C(\alpha) \frac{|\nabla \eta|^4}{\eta^3}f^{2\beta+4} \nonumber \\
&\le \frac{(\alpha -1)}{4} \eta F^2 + \frac{C(\alpha)}{R^4} f^{2\beta+4}, 
\end{align}
where we have used $\sqrt F = f^{\beta/2} |\nabla f|$ and in much the same way 
\begin{align}
2 \frac{|\nabla \eta|^2}{\eta}f^{\beta+2} F = 2 \eta^{1/2} F \frac{|\nabla \eta|^2}{\eta^{3/2}} f^{\beta+2} 
&\le \frac{\alpha -1}{4} \eta F^2 + C(\alpha) \frac{|\nabla \eta|^4}{\eta^3} f^{2\beta+4} \nonumber\\
&\le \frac{\alpha -1}{4} \eta F^2 + \frac{C(\alpha)}{R^4} f^{2\beta+4}. 
\end{align}
\item For the terms involving $\Delta_\phi \eta$ and $\partial_t \eta$ as in the proof of Theorem \ref{thm1} 
we have the bounds  
\begin{align*}
- \Delta_\phi \eta &  \leq C \left( \frac{1}{R^2} + \frac{[\gamma_{\Delta_\phi}]_+}{R} + (n-1)k \right) \sqrt \eta \qquad 
\text{and } \qquad \partial_t \eta \le  \frac{c\sqrt \eta}{\tau-t_0+T}.
\end{align*}
Hence, by substituting and making use of Young's inequality we have, 
\begin{align} \label{Lap comp 2 sum}
&-f^{\beta+2}F \left( \Delta_\phi - \partial_t -2(n-1)k \right) \eta \\
&\le C(n) F \sqrt \eta 
\left( \frac{1}{\tau-t_0+T} + k + \frac{1}{R^2} 
+ \frac{[\gamma_{\Delta_\phi}]_+}{R}
 \right) f^{\beta+2} \nonumber \\
& \leq \frac{\alpha-1}{4} \eta F^2
+ C(\alpha, n) \left( 
\frac{1}{(\tau-t_0+T)^2} + k^2 
+ \frac{1}{R^4} + \frac{[\gamma_{\Delta_\phi}]_+^2}{R^2}
\right) f^{2\beta+4}. \nonumber
\end{align}
\item Making note of 
$|\langle \nabla f, {\mathscr G}_x(t,x, f^\alpha) \rangle| \le |\nabla f| |\mathscr G_x(t,x, f^\alpha)|$ 
and $w=f^{\alpha}$ we have by abbreviating the arguments $(t,x, f^\alpha)$
\begin{align}\label{eq.5.18}
&\frac{2w \mathscr G_w - [2-(2+\beta)/\alpha] \mathscr G}{w} f^{\beta+2} \eta F
+ \frac{(2/\alpha) f^{2\beta+3} \langle \nabla f, \mathscr G_x \rangle \eta}{w} \nonumber \\
&\le \left[ \frac{2w \mathscr G_w - [2-(2+\beta)/\alpha] \mathscr G}{w} \right]_+ f^{\beta+2} \eta F
+ (2/\alpha) f^{3\beta/2+3} \sqrt F \frac{|\mathscr G_x|}{w} \eta \nonumber \\
&\le \frac{\alpha -1}{4} \eta F^2 + C(\alpha) f^{2\beta+4} \left\{ \left [\frac{|\mathscr G_x(t,x,w)|}{w}\right]^{4/3} \right. \nonumber\\
&\left. \,\,\,\, + \left[\frac{2w \mathscr G_w(t,x,w)-[2-(2+\beta)/\alpha] \mathscr G(t,x,w)}{w} \right]_+^2 \right\}, 
\end{align}
\end{itemize}

Now putting together \eqref{eq.5.15}-\eqref{eq.5.18} it follows from \eqref{eq4.2} 
after some basic calculations and rearrangement 
of terms that at the space-time point $(x_1, t_1)$ we have

\begin{align} \label{ing eq1}
[\alpha -1 -\beta^2 &-\beta(2-\alpha)] \eta F^2 
\le~ C(\alpha, n) \Big( \sup_{Q_{R, T}}  w \Big)^{(2\beta+4)/\alpha} \times \nonumber\\
&\left \{ \begin {array}{ll} k^2+\dfrac{1}{(\tau-t_0+T)^2} + \dfrac{1}{R^4} \\
+ \dfrac{[\gamma_{\Delta_\phi}]_+^2}{R^2} 
+ \sup_{Q_{R, T}}  \left [\dfrac{|\mathscr G_x(t,x,w)|}{w}\right]^{4/3} \\
+ \sup_{Q_{R, T}}  \left[\dfrac{2w \mathscr G_w(t,x,w)-[2-(2+\beta)/\alpha] \mathscr G(t,x,w)}{w} \right]_+^2
\end{array}
\right\}.
\end{align}

By the maximality of $\eta F$ at $(x_1, t_1)$ we have for any $x$ with $d(x, x_0) \le R/2$ the chain of inequalities 
$F(x,\tau) \le [\eta F](x, \tau) \le [\eta F](x_1,t_1) \le [\sqrt\eta F](x_1,t_1)$ [recall that $\eta(x, \tau) =1$ 
when $d(x,x_0) \leq R/2$]. Hence combining the latter with \eqref{ing eq1} and making note of the relations 
$f^\alpha=w$ and $F = f^\beta |\nabla f|^2 =  w^{(\beta+2)/\alpha} |\nabla w|^2/(\alpha ^2w^2)$ this gives 
\begin{align}
s(\alpha, \beta)\frac{|\nabla w|^2}{w^{2-(\beta+2)/\alpha}} =&~ s(\alpha, \beta) \alpha^2 f^\beta |\nabla f|^2
\nonumber\\
\le&~ C(\alpha, n) \Big( \sup_{Q_{R, T}}  w \Big)^{(\beta+2)/\alpha} \times \nonumber \\
&\left \{ \begin {array}{ll} 
k+ \dfrac{1}{\tau-t_0+T} + \dfrac{1}{R^2} \\
+ \dfrac{[\gamma_{\Delta_\phi}]_+}{R} 
+ \sup_{Q_{R, T}}  \left [\dfrac{|\mathscr G_x(t,x,w)|}{w}\right]^{2/3}\\
+ \sup_{Q_{R, T}}  \left[\dfrac{2w \mathscr G_w(t,x,w)-[2-(2+\beta)/\alpha] \mathscr G(t,x,w)}{w} \right]_+
\end{array}
\right\},
\end{align}
where we have set
\begin{align}
s(\alpha, \beta) = \frac{[\alpha -1 -\beta^2 -\beta(2-\alpha)]^{1/2}}{\alpha^2} &= \frac{[1 -(1 +\beta^2)/\alpha -(2/\alpha-1)\beta]^{1/2}}{\alpha^{3/2}}\nonumber\\
& =\frac{[(1+\beta)(1-(1+\beta)/\alpha)]^{1/2}}{\alpha^{3/2}}.
\end{align}
Note that the condition $0<(1+\beta)<\alpha$ guarantees that $(1+ \beta)[1-(1+\beta)/\alpha]>0$ 
and so in turn $s(\alpha, \beta)>0$. This upon rearranging terms, taking square roots and noting the 
arbitrariness of $\tau$ gives the desired estimate for every $(x, t) \in{\mathscr{B}_{R/2}\times (t_0-T, t_0]}$. \qed

\section{Proof of the Liouville Theorems and implications $({\bf I})$: ${\mathscr Ric}_\phi(g) \ge 0$} \label{sec6}

This section is devoted to the proofs of the Liouville constancy results formulated in Theorem \ref{Liouville log thm} 
through to Theorem \ref{ancient}. Note that in Theorem \ref {Liouville log thm} to Theorem \ref{LiouvilleThmEx-2.9} 
it is the elliptic counterpart of $({\mathscr P})$ with ${\mathscr G}={\mathscr G}(w)$ that is being considered and 
so evidently here the solution $w$ is time independent.

\qquad \\
{\bf Proof of Theorem \ref{Liouville log thm}.} Let $h=\log(w/D)$ and note that in view of $0<w \le D$ we have $h \le 0$ 
and so $1/(1-h) \le 1$. Referring to the estimate \eqref{eq13}, making note of ${\mathscr G}={\mathscr G}(w)$, the solution 
$w$ being time independent, $k=0$ and  then passing to the limit $t \nearrow \infty$ (hence $1/\sqrt{t-t_0+T} \to 0$) gives 
for each fixed $x$ in ${\mathscr B}_{R/2}$ the bound 
\begin{align} 
\frac{|\nabla w|}{w}(x)   
\le &~C  \left(1 - \log \frac{w}{D}(x) \right) \times \\
& \left \{ \begin {array}{ll}  \sqrt{k} + \dfrac{1}{R} + \sqrt{\dfrac{[\gamma_{\Delta_\phi}]_+}{R}} 
+ \sup_{{\mathscr B}_{R}} \left[\dfrac{|\mathscr G_x(w)|}{ w[1- \log(w/D)]^2}\right]^\frac{1}{3} \\
+  \sup_{{\mathscr B}_{R}} \left[ \dfrac{w[1-\log(w/D)] \mathscr G_w(w) 
+ \log(w/D) \mathscr G (w)}{w [1-\log(w/D)]^2} \right]_+^\frac{1}{2} \nonumber 
\end{array} 
 \right\} \\
 \le &~C (1-h(x)) [1/R + [\gamma_{\Delta_\phi}]_+^{1/2}/\sqrt R 
 + \sup_{{\mathscr B}_R} \{[(1-h) \mathscr G'(w) + h \mathscr G (w)/w]_+^{1/2} \}]. \nonumber 
\end{align}
From $(1-h) w \mathscr G'(w) + h \mathscr G (w) \le 0$ it follows that 
$[(1-h) \mathscr G'(w) + h \mathscr G (w)/w]_+ \equiv 0$ and passing to the limit 
$R \nearrow \infty$ gives $|\nabla w|(x) =0$. The arbitrariness of $x\in M$ now gives 
$|\nabla w| \equiv 0$ on $M$ and the connectedness of $M$ implies that $w$ is a constant. 
\hfill $\square$

\qquad \\
{\bf Proof of Theorem \ref{Liouville two thm}.}
As the methodology and main ideas here are to some extenet similar to those used in the proof 
of Theorem \ref{Liouville log thm}, we shall remain brief, focusing mainly on the 
differences.  Referring to the estimate \eqref{eq1.12}, making note of ${\mathscr G}={\mathscr G}(w)$, the solution 
$w$ being time independent, $k=0$ and  then passing to the limit $t \nearrow \infty$ (hence $1/\sqrt{t-t_0+T} \to 0$) 
we have for each fixed $x$ in ${\mathscr B}_{R/2}$ the bound
\begin{align}
\frac{|\nabla w|}{w^{1-(\beta+2)/(2\alpha)}}(x)
\le&~C \Big( \sup_{{\mathscr B}_{R}} w \Big)^{(\beta+2)/(2\alpha)} \times \\
&\left \{ \begin {array}{ll} 
\sqrt k + \dfrac{1}{R} + \sqrt{\dfrac{[\gamma_{\Delta_\phi}]_+}{R}}
+\sup_{{\mathscr B}_{R}} \left[|\mathscr G_x (w)|/w\right]^\frac{1}{3} \nonumber \\
+\sup_{{\mathscr B}_{R}} \left[2\mathscr G_w(w) -
[2-(\beta+2)/\alpha] \mathscr G(w)/w \right]_+^\frac{1}{2}
\end {array}
\right\}, \nonumber \\
&\le C \Big( \sup_{{\mathscr B}_{R}} w \Big)^{(\beta+2)/(2\alpha)} 
[1/R + [\gamma_{\Delta_\phi}]_+^{1/2}/\sqrt R \nonumber \\
&+ \sup_{{\mathscr B}_{R}} \{ [2\mathscr G'(w) -
[2-(\beta+2)/\alpha] \mathscr G(w)/w]_+^\frac{1}{2} \}]. \nonumber
\end{align}
Now $[1-(\beta/2+1)/\alpha] {\mathscr G}(w) -w {\mathscr G}'(w) \ge 0$ gives 
$[2\mathscr G'(w) - [2-(\beta+2)/\alpha] \mathscr G(w)/w]_+ \equiv 0$
and making note of $\sup w< \infty$ and $(\beta+2)/(2\alpha)>0$ and passing to the limit 
$R \nearrow \infty$ gives $|\nabla w|(x) =0$. 
The arbitrariness of $x\in M$ now gives $|\nabla w| \equiv 0$ on $M$ and again the connectedness of $M$ 
implies that $w$ is a constant. 
\hfill $\square$

\qquad \\
{\bf Proof of Theorem \ref{Liouville one thm}.}
Consider the function ${\mathscr G}(w) =  X(w) + w^r Y(\log w)$ where $w>0$. A basic differentiation then gives 
${\mathscr G}'(w) = X'(w) + w^{r-1} [rY(\log w) + Y'(\log w)]$ and by so using the sub-additivity of $[\cdot]_+$ and recalling 
that $h=\log (w/D)$ we can write
\begin{align} \label{app one R}
\mathsf L =&~ \left[ \frac{w (1-h) {\mathscr G}'(w) + h {\mathscr G}(w)}{(1-h)^2 w} \right]_+ 
\le \left[ \frac{\{ w (1-h) X'(w) + h X(w) \}}{(1-h)^2 w} \right]_+ \\
&+ \left[ \frac{ w^r \{ [(1-h) r +h] Y(\log w) + (1-h) Y'(\log w) \}}{(1-h)^2 w} \right]_+ 
= \mathsf{L}_1+ \mathsf{L}_2.  \nonumber 
\end{align}

By inspection it is seen that $(1-h) w X' + h X = w X' - h (wX' - X) \le 0$
for $w>0$ as a result of the assumptions on $X$ in the theorem.
Therefore
\begin{align} \label{Psi term []+}
\mathsf{L}_1 &= \left[ \frac{ (1-h) w X'(w) + h X (w)}{(1-h)^2 w} \right]_+ 
= \left[ \frac{w X'(w) + h [X (w) - w X'(w)]}{(1-h)^2 w} \right]_+ \equiv 0.   
\end{align}  
Next we have $[r+h/(1-h)] Y(\log w) + Y'(\log w) \le 0$ for $w \ge 1$, $r \ge 1$ 
when $Y$ satisfies $({\mathsf Y1})$. Thus here we have 
\begin{align}
\mathsf{L}_2
&= \left[ \frac{ w^r \{ [(1-h) r +h] Y(\log w) + (1-h) Y'(\log w) \}}{(1-h)^2 w} \right]_+ \nonumber \\
&= \frac{w^{r-1}}{1-h}  \left[ \left( r + \frac{h}{1-h} \right) Y(\log w) + Y'(\log w) \right]_+ \equiv 0. 
\end{align}
Similarly from $({\mathsf Y2})$ we have $s Y'(s) \ge \gamma Y(s)$ for $s\le 0$ and so when 
$r \le \min (\gamma, 1)$ and $0<w \le 1$ we have 
\begin{align} 
\frac{r - h(r-1)}{1-h} & Y(\log w) + Y'(\log w) = \frac{[r - h(r-1)] Y(\log w) + (1-h) Y'(\log w)}{1-h} \nonumber \\
&\le [r -\gamma - h(r-1)] \frac{Y(\log w)}{1-h} + (1+\log D) \frac{Y'(\log w)}{1-h} \le 0, 
\end{align}
by writing $h=\log(w/D)$ and choosing $D \ge 1$. Therefore again here we have 
\begin{align}
\mathsf{L}_2
&= \left[ \frac{ w^r \{ [(1-h) r +h] Y(\log w) + (1-h) Y'(\log w) \}}{(1-h)^2 w} \right]_+ \nonumber \\
&= \frac{w^{r-1}}{1-h}  \left[\frac{r - h(r-1)}{1-h} Y(\log w) + Y'(\log w) \right]_+ \equiv 0.
\end{align} 

Thus referring to in \eqref{app one R} we have ${\mathsf L}= 0$. 
An application of the estimate \eqref{eq13} as in the proof of Theorem \ref{Liouville two thm} now 
gives the desired conclusion. 
\hfill $\square$

\qquad \\
{\bf Proof of Theorem \ref{Liouville.2-8}.}
The idea is to invoke Theorem \ref{Liouville two thm} with the choices $X=X(w)$ as in \eqref{polies PR} and $Y \equiv 0$. 
Indeed referring to the formulation of Theorem \ref{Liouville two thm} it is seen by a direct calculation that 
\begin{equation} 
X'(w) = \sum_{j=1}^N p_j \mathsf{A}_j w^{p_j-1} + \sum_{j=1}^N q_j \mathsf{B}_j w^{q_j-1},  
\end{equation}
and thus 
\begin{align}
[1-(\beta/2+1)/\alpha]X(w) -w X'(w) 
=&~\sum_{j=1}^N \mathsf{A}_j [1-(\beta/2+1)/\alpha - p_j] w^{p_j} \nonumber \\
&+ \sum_{j=1}^N \mathsf{B}_j [1-(\beta/2+1)/\alpha - q_j] w^{q_j}. 
\end{align}
This is then easily seen to be non-negative, as required by Theorem  \ref{Liouville two thm}, upon suitably restricting the 
ranges of $\mathsf{A}_j$, $\mathsf{B}_j$ and $p_j, q_j$ as formulated in the statement of the theorem.

\begin{remark}
Note that the upper and lower bounds on the exponents $p_j, q_j$ are given by the same quantity 
$1-(\beta/2+1)/\alpha$ which can be adjusted by optimising the parameters $\alpha, \beta$ within 
their respective range. In particular if $\sum {\mathsf B}_j w^{q_j} \equiv 0$ we only require $\mathsf{A}_j \ge 0$ 
and $p_j <1$ $($with $\alpha=0$, $\beta \searrow 0$$)$ and if $\sum {\mathsf A}_j w^{p_j} \equiv 0$ we only require 
$\mathsf{B}_j \le 0$ and $q_j >0$ $($with $\alpha=0$, $\beta \nearrow 1$$)$. 
\end{remark}

\qquad \\
{\bf Proof of Theorem \ref{LiouvilleThmEx-2.9}.}
Writing ${\mathscr G}(w) = \mathsf{A} w^p + \mathsf{B} w^q + w^r Y(\log w)$ a straightforward calculation gives
\begin{align}
[1-(\beta/2+1)/\alpha] \mathscr G (w) - w \mathscr G' (w)
=&~ w^r ([1-(\beta/2+1)/\alpha -r] Y - Y') \nonumber\\
&+ \mathsf{A} [1-(\beta/2+1)/\alpha-p] w^p \nonumber\\
&+ \mathsf{B} [1-(\beta/2+1)/\alpha - q] w^{q}.
\end{align}  
The conclusion follows from an application of the first part of Theorem \ref{Liouville two thm} upon noting that according 
to assumptions $[1-(\beta/2+1)/\alpha] \mathscr G (w) - w \mathscr G' (w) \ge 0$. \hfill $\square$

\qquad \\
{\bf Proof of Theorem \ref{ancient}.}
Fix a space-time point $(x_0, t_0)$. Then with the quantities $k=0$, $t=t_0$ and $T=R$ in \eqref{eq13} and upon writing 
$D_{R,R} = \sup_{Q_{R,T}} w$ and making note of $\log w(x,t) = o[\sqrt{r(x)} + \sqrt{|t|}]$ we have
\begin{align} 
\frac{|\nabla w|}{w} (x_0,t_0) 
\le &~C  \left(1 - \log w (x_0,t_0) + \log D_{R,R} \right) \times \\
& \left \{ \begin {array}{ll}  \sqrt{k}+\dfrac{1}{\sqrt{t-t_0+T}} + \dfrac{1}{R}\\
+\sqrt{\dfrac{[\gamma_{\Delta_\phi}]_+}{R}}
+ \sup_{Q_{R, T}} \left[\dfrac{|\mathscr G_x(w)|}{w[1- \log(w/D)]^2}\right]^\frac{1}{3} \\
+  \sup_{Q_{R, T}} \left[ \dfrac{w[1-\log(w/D)] \mathscr G_w(w) 
+ \log(w/D) \mathscr G (w)}{w [1-\log(w/D)]^2} \right]_+^\frac{1}{2} \nonumber 
\end{array} 
\right\} \nonumber \\
\le &~C  \left(1 - \log w (x_0,t_0) + \log D_{R,R} \right)[1/\sqrt R + 1/R + [\gamma_{\Delta_\phi}]^{1/2}_+/\sqrt R] \nonumber 
\end{align}
where as in the proof of Theorem \ref{Liouville one thm} we have deduced from the assumptions ${\mathscr G}' (w) \le 0$ 
and ${\mathscr G}(w) - w {\mathscr G}'(w) \ge 0$ that 
$[  (1-\log(w/D)) w \mathscr G'(w) + \log(w/D) \mathscr G (w)]_+=0$. Passing to the limit $R \nearrow \infty$ and using the 
assumption on the growth of $w=w(x,t)$ it follows that $|\nabla w(x_0,t_0)|=0$. The arbitrariness of $(x_0,t_0)$ implies 
$|\nabla w| \equiv 0$ and so $w=w(t)$. From the equation \eqref{ancient-equation} it then follows that 
$dw/dt = {\mathscr G}(w)$. Integrating the latter and using the assumption ${\mathscr G}(w) \ge a >0$ for all $w>0$ gives $w(t) \le w(0) + at$ 
for all $t<0$. This however clashes with $w(t)>0$ as $t \searrow -\infty$ and so the conclusion is reached. \hfill $\square$

\section{Proof of the local Li-Yau estimate in Theorem \ref{thm7.1}} \label{sec7}

We now attend to the proof of the local Li-Yau estimate (also called a differential Harnack 
inequality) as formulated in Theorem \ref{thm7.1}. 
The first two subsections develop the necessary tools including an evolution inequality on 
a suitably defined Harnack quantity (below called $F$) built out of the solution $w$, and involving a parameter $\alpha>1$, and the last 
section completes the proof by putting these together, using a localisation argument and eventually maximum principle.

\subsection{Evolution of the Harnack quantity $F=  t [|\nabla w|^2/w^2 - \alpha \partial_t w/w + \alpha \mathscr G/w]$}

From the positive solution $w$ we first define the function $f=\log w$ and then define the Harnack quantity 
$F= t [|\nabla f|^2 - \alpha \partial_t f + \alpha e^{-f} \mathscr G(t,x,e^f)]$ (for all $t \ge 0$ and fixed non-zero 
parameter $\alpha$). The evolution of $F=F_\alpha$ under the equation $({\mathscr P})$ is the subject of the next lemma.

\begin{lemma}  \label {estimate-second-type-lemma} 
Let $w$ be a positive solution to $({\mathscr P})$. For $f=\log w$ and $\alpha$ a fixed non-zero constant let 
$F=F_\alpha(x,t)$ be defined by 
\begin {align} \label{2.3}
F &=  t \left[ \frac{|\nabla w|^2}{w^2} - \alpha \frac{\partial_t w}{w} + \alpha \frac{\mathscr G(t,x,w)}{w} \right] 
=   t \left[ \frac{|\nabla w|^2}{w^2} - \alpha \frac{\Delta_\phi w}{w} \right]  \nonumber \\
&= t [|\nabla f|^2 - \alpha \partial_t f + \alpha e^{-f} \mathscr G(t,x,e^f)], \qquad t \ge 0. 
\end{align}  
Then $F$ satisfies the evolution equation 
\begin{align} \label{9.4a}
(\Delta_\phi - \partial_t)F =& ~2 t|\nabla \nabla f|^2 - 2 \langle \nabla f , \nabla F \rangle \nonumber\\
&+2t {\mathscr Ric}_\phi^m (\nabla f, \nabla f) + [2t/(m-n)] \langle \nabla \phi , \nabla f \rangle^2 \nonumber \\
& -F/t+2t (\alpha -1) \langle \nabla f, \nabla [e^{-f}\mathscr G(t,x,e^f)] \rangle \nonumber\\
&+ \alpha t \Delta_\phi [e^{-f} \mathscr G(t,x,e^f)].
\end{align}
In particular, if ${\mathscr Ric}_\phi^m(g) \ge -(m-1) kg$ for some $n \le m <\infty$ and $k \ge 0$ then 
\begin {align} \label{eq7.9}
(\Delta_\phi- \partial_t)F \ge&~2t(\Delta_\phi f)^2/m 
- 2 \langle \nabla f, \nabla F \rangle \nonumber \\
&- 2 t (m-1)k |\nabla f|^2 - F/t \nonumber\\
&+2 t (\alpha -1)\langle \nabla f, \nabla (e^{-f}\mathscr G(t,x,e^f)) \rangle \nonumber \\
&+\alpha t \Delta_\phi [e^{-f}\mathscr G(t,x,e^f)].
\end{align}
\end{lemma}

\begin{proof}
From $({\mathscr P})$ it follows that $f=\log w$ satisfies the equation 
\begin {align} \label{2.2}
\partial_t f= \Delta_\phi f + |\nabla f|^2 + e^{-f} \mathscr G(t,x,e^f).
\end{align} 
Thus using \eqref{2.3} and \eqref{2.2} it is seen that 
\begin {align} \label{9.6a}
\Delta_\phi f 
&= \partial_t f - |\nabla f|^2 - e^{-f} \mathscr G(t,x,e^f) \nonumber \\
&= - \left[ (1/\alpha) |\nabla f|^2 - \partial_t f + e^{-f} \mathscr G(t,x,e^f) \right] 
- [(\alpha -1)/\alpha] |\nabla f|^2 \nonumber \\
&= -F/(\alpha t) - [(\alpha-1)/\alpha] |\nabla f|^2, \qquad t>0. 
\end{align}
Next, an application of $\Delta_\phi$ to $F$ in \eqref{2.3} gives   
\begin {align} \label{Witten-Lap-F-DHI-equation}
\Delta_\phi F =&~t (\Delta_\phi |\nabla f|^2 - \alpha \Delta_\phi (\partial_t f) + \alpha \Delta_\phi [e^{-f}\mathscr G(t,x,e^f)]) \nonumber \\
=&~2t |\nabla \nabla f|^2 +2t \langle \nabla f, \nabla \Delta_\phi f \rangle 
+2t {\mathscr Ric}_\phi^m (\nabla f, \nabla f) \nonumber\\
&+[2t/(m-n)] \langle \nabla \phi , \nabla f \rangle^2 - \alpha t \partial_t (\Delta_\phi f) 
+ \alpha t \Delta_\phi [e^{-f} \mathscr G (t,x,e^f)]
\end{align} 
where in concluding the second line we have used the weighted Bocnher-Weitzenb\"ock formula as applied to $f$. 
Hence by virtue of the  \eqref{2.2} and \eqref{9.6a} this can be simplified and re-written as 
\begin {align} \label{2.11}
\Delta_\phi F 
=&~ 2t |\nabla \nabla f|^2 +2t {\mathscr Ric}_\phi^m (\nabla f, \nabla f)
+[2t/(m-n)] \langle \nabla \phi , \nabla f \rangle^2
+2t \langle \nabla f, \nabla \Delta_\phi f \rangle \nonumber\\
&+(t\partial_t F - F)/t +2t (\alpha -1)\langle \nabla f , \nabla (\partial_t f)\rangle 
+ \alpha t \Delta_\phi [e^{-f} \mathscr G (t,x,e^f)]\nonumber\\
=& ~ 2t |\nabla \nabla f|^2 +2t {\mathscr Ric}_\phi^m (\nabla f, \nabla f)
+[2t/(m-n)] \langle \nabla \phi , \nabla f \rangle^2 \nonumber\\
& +2t (\alpha -1) \langle \nabla f , \nabla [e^{-f}\mathscr G (t,x,e^f)] \rangle + \partial_t F - F/t \nonumber\\
&- 2 \langle \nabla f , \nabla F \rangle + \alpha t \Delta_\phi [e^{-f} \mathscr G (t,x,e^f)]. 
\end{align}

Rearranging this now lead to 
\begin{align}\label{eq3.8}
(\Delta_\phi - \partial_t) F = &~2 t|\nabla \nabla f|^2  + 2t {\mathscr Ric}_\phi^m (\nabla f, \nabla f)\nonumber\\
&+ [2t/(m-n)] \langle \nabla \phi , \nabla f \rangle^2 - F/t -2 \langle \nabla f, \nabla F \rangle \nonumber\\
&+ 2t (\alpha -1) \langle \nabla f, \nabla[e^{-f}\mathscr G(t,x,e^f)] \rangle
+ \alpha t \Delta_\phi [ e^{-f} \mathscr G (t,x,e^f)],
\end{align}
which is the desired conclusion.

For the inequality in  \eqref{eq7.9} we first note that by an application of the Cauchy-Schwarz 
and Young inequalities we have 
\begin{equation}
|\nabla \nabla f|^2+\langle \nabla \phi , \nabla f \rangle^2/(m-n) 
\ge (\Delta f)^2/n+\langle \nabla \phi , \nabla f \rangle^2/(m-n) 
\ge (\Delta_\phi f)^2/m.
\end{equation}
Hence the desired conclusion follows by substituting the above in \eqref{eq3.8} and applying the Ricci curvature lower bound 
${\mathscr Ric}_\phi^m(g) \ge -(m-1) k g$ as given in the lemma.
\end{proof}

In the course of the proof of Theorem \ref{thm7.1} we make use of the following lemma. 
As before we use the notation $\mathscr G^x: x \mapsto \mathscr G(t,x, w)$ for the function obtained by freezing 
the variables $t$, $w$ and viewing $\mathscr G$ solely as a function of $x$. We also denote partial derivatives of 
${\mathscr G}$ with subscripts.

\begin{lemma} 
Suppose $\mathscr G=\mathscr G(t,x,w)$ is a twice continuously differentiable function on 
$[0, T] \times M \times (0, \infty)$ and let $w=e^f$ with $f=f(x,t)$ be a twice continuously 
differentiable on $M \times [0, T]$. Then 
\begin{align} 
\Delta_\phi \mathscr G 
=  \Delta_\phi \mathscr G^x +2e^f \langle \mathscr G_{xw} ,\nabla f \rangle 
+ e^f |\nabla f|^2 ( \mathscr G_w +e^f \mathscr G_{ww})+e^f \mathscr G_w \Delta_\phi f.
\end{align}
\end{lemma}

\begin{proof}
Firstly, it is seen, by calculating in local coordinates or directly otherwise that with 
${\mathscr G}={\mathscr G}(t,x,w)$ and $w=e^f$ we have 
\begin {align} \label{3.31}
\nabla \mathscr G = \mathscr G_x + e^f \mathscr G_w \nabla f, 
\qquad \mathscr G_x = (\mathscr G_{x_1}, \dots, \mathscr G_{x_n}), 
\end{align}  
where we have abbreviated the arguments $(t, x, w)$ for convenience. Now differentiating further 
with the aim of calculating the Laplacian (and noting $\mathscr G_x= \nabla \mathscr G^x$) gives 
\begin{align}
\Delta \mathscr G 
&= \Delta \mathscr G^x + e^f \langle \mathscr G_{xw} ,\nabla f \rangle +e^f |\nabla f|^2 \mathscr G_w 
+ e^f \langle \mathscr G_{xw},\nabla f \rangle
+e^{2f} |\nabla f|^2 \mathscr G_{ww} +e^f \mathscr G_w \Delta f \nonumber\\
& = \Delta \mathscr G^x + 2e^f \langle \mathscr G_{xw} ,\nabla f \rangle 
+ e^f |\nabla f|^2 (\mathscr G_w +e^f \mathscr G_{ww}) + e^f \mathscr G_w \Delta f.
\end{align}   
Subsequently calculating the Witten Laplacian, by utilising the above fragments we have,  
\begin {align} \label{3.33}
\Delta_\phi \mathscr G &= \Delta \mathscr G -\langle \nabla \phi , \nabla \mathscr G \rangle = \Delta \mathscr G 
-\langle \nabla \phi , ( \mathscr G_x+e^f \mathscr G_w\nabla f) \rangle \\
& = \Delta \mathscr G -\langle \nabla \phi, \mathscr G_x \rangle -e^f \mathscr G_w \langle \nabla \phi, \nabla f \rangle\nonumber \\
& = \Delta \mathscr G^x -\langle \nabla \phi, \nabla \mathscr G^x \rangle +2e^f \langle \mathscr G_{xw} ,\nabla f \rangle 
+ e^f |\nabla f|^2 ( \mathscr G_w +e^f \mathscr G_{ww})+e^f \mathscr G_w \Delta_\phi f, \nonumber
\end{align}
which is the required result.
\end{proof}

\subsection{Construction of spatial cut-offs and localisation in space}

The next and final ingredient needed in the proof of Theorem \ref{thm7.1} is the construction of a suitable cut-off function in space. 
To this end we pick a reference point $x_0 \in M$ and fix $R>0$ and with $r(x)=d(x,x_0)$ being the 
geodesic radial variable with respect to $x_0$ we set 
\begin{align}\label{9.31}
\zeta(x)=\bar{\zeta} \left( \frac{r (x)}{R} \right). 
\end{align}
The function $\bar{\zeta}=\bar \zeta(s)$ appearing in \eqref{9.31} is defined on the half-line $s \ge 0$ and can be 
constructed easily with properties described below (see \cite{[LY86], [LiP], [Wu18]}).

\begin{lemma} \label{psi lemma} There exists a function $\bar{\zeta}:[0,\infty) \to \mathbb{R}$ such that:
\begin{enumerate}[label=$(\roman*)$]
\item $\bar\zeta$ is of class $\mathscr{C}^2 [0, \infty)$, 
\item $0 \le \bar\zeta(s) \le 1$ for $0 \le s < \infty$ with $\bar{\zeta}(s) \equiv 1$ for $s \le 1$ and $\bar{\zeta}(s) \equiv 0$ for $s \ge 2$.
\item $\bar\zeta$ is non-increasing, $($$\bar\zeta' \le 0$$)$ and additionally for suitable constants $c_1, c_2>0$, 
\begin{align} \label{9.30}
- c_1 \le \frac{\bar{\zeta}^{'}}{\sqrt{ \bar{\zeta}}} \le 0 \qquad and \qquad  \bar{\zeta}^{''} \ge -c_2. 
\end{align}
\end{enumerate}
\end{lemma}

It is evident that here $\zeta \equiv 1$ for when $0 \le r(x) \le R$ and $\zeta \equiv 0$ for when $r(x) \ge 2R$. 
Furthermore a straightforward calculation starting from \eqref{9.31} gives 
\begin{align} \label{eq3.25}
\nabla \zeta = (\bar\zeta'/R) \nabla r, \qquad
\Delta \zeta = \bar\zeta'' |\nabla r|^2/R^2+\bar\zeta' \Delta r/R
\end{align}
and so 
\begin{align} \label{eq3.25-alt}
[\Delta_\phi-\partial_t] \zeta = \Delta_\phi \zeta = \bar\zeta'' |\nabla r|^2/R^2+\bar\zeta' \Delta_\phi r/R.
\end{align}

\subsection{Finalising the proof of Theorem \ref{thm7.1}}

We now come to the proof of the main result in the section which is the local Li-Yau estimate in Theorem \ref{thm7.1}. To this end 
we consider the spatially localised function $\zeta F$ where $F$ is as in \eqref{2.3}. 
We denote by $(x_1,t_1)$ the point where this function attains its maximum over the compact 
cylinder $\{d(x,x_0) \le 2R, 0 \le t \le T\}$. We also assume that $[\zeta F] (x_1, t_1)>0$ as otherwise 
the desired estimate is trivially true as a result of $F \le 0$. It thus follows that $t_1>0$ and 
$d(x_1, x_0) < 2R$ and so at the maximum point $(x_1, t_1)$ we have the relations 
\begin{align} \label{9.32}
\begin{cases}
\partial_t (\zeta F) \ge 0, \\
\nabla (\zeta F)=0, \\
\Delta (\zeta F) \le 0.
\end{cases}
\end{align}
From the inequalities above we deduce that  
\begin{equation}
0 \ge [\Delta_\phi -\partial_t] (\zeta F) = F [\Delta_\phi -\partial_t] \zeta + 2 \langle \nabla \zeta, \nabla F \rangle + \zeta [\Delta_\phi - \partial_t] F.
\end{equation}

Now making note of $[\Delta_\phi -\partial_t] \zeta = \Delta_\phi \zeta$ (as $\partial_t \zeta \equiv 0$) it follows from the above that at the 
maximum point $(x_1, t_1)$ we have 
\begin{align} \label{9.34}
0 &\ge F \Delta_\phi \zeta + 2 \langle \nabla \zeta, \nabla F \rangle + \zeta [\Delta_\phi - \partial_t] H\nonumber \\
& \ge F \Delta_\phi \zeta + (2/\zeta) \langle \nabla \zeta , \nabla (\zeta F) \rangle -2 (|\nabla \zeta |^2/\zeta) F + \zeta [\Delta_\phi - \partial_t ]F\nonumber \\
& \ge F \Delta_\phi \zeta -2 (|\nabla \zeta |^2/\zeta) F + \zeta [\Delta_\phi - \partial_t] F.
\end{align}

Next by referring to \eqref{eq3.25}-\eqref{eq3.25-alt} and $\Delta_\phi r \le (m-1) \sqrt {k}\coth (\sqrt{k} r)$
[the latter being a consequence of the generalised Laplacian comparison theorem and the lower curvature 
bound ${\mathscr Ric}_\phi^m \geq -(m-1) k g$] it follows upon recalling $\bar\zeta'\le0$ in $(iii)$ in Lemma \ref{psi lemma} that 
\begin{align}
\Delta_\phi \zeta 
= \bar\zeta'' \frac{|\nabla r|^2}{R^2}+\bar\zeta' \frac{\Delta_\phi r}{R}
\ge \frac{1}{R^2} \bar\zeta''+\frac {(m-1)}{R} \bar\zeta' \sqrt{k} \coth (\sqrt{k} r).
\end{align}

Moreover upon noting $\coth (\sqrt{k} r) \le \coth (\sqrt{k} R)$ and $\sqrt{k} \coth (\sqrt{k} R) \le (1 + \sqrt{k} R)/R$, 
for $R \le r \le 2R$, we deduce that $(m-1) \bar\zeta' \sqrt{k} \coth (\sqrt{k}r)\ge (m-1)(1+\sqrt{k}R)/R \bar\zeta$. 
Hence by putting the above fragments together we have the lower bound on $\Delta_\phi \zeta$ in terms of $k$, $R$ 
and the constants $c_1$, $c_2$ in \eqref{9.30}: 
\begin{align} \label{9.38}
\Delta_\phi \zeta 
&\ge \frac{1}{R^2} \bar\zeta''+\frac {(m-1)}{R} \bar\zeta' \sqrt{k} \coth (\sqrt{k} r) 
\ge \frac {1}{R^2}\bar{\zeta}^{''}+\frac{(m-1)}{R} \left( \frac{1}{R}+\sqrt{k} \right) \bar{\zeta}^{'}\nonumber\\
&\ge -\frac{c_2}{R^2}-\frac{(m-1)}{R} c_1 \left(\frac{1}{R}+\sqrt{k}\right) 
=  -\frac{1}{R^2}[c_2 +(m-1) c_1(1+R \sqrt{k})].
\end{align}
Likewise, referring to \eqref{9.30} and \eqref{eq3.25}, a straightforward calculation gives,
\begin{align} \label{9.39}
\frac{|\nabla \zeta |^2}{\zeta} = \frac{\bar\zeta'^2}{\bar\zeta} \frac{|\nabla r|^2}{R^2} 
= \left( \frac{\bar\zeta'}{\sqrt {\bar\zeta}} \right)^2 \frac{|\nabla r|^2}{R^2} \le  \frac{c_1^2}{R^2}.  
\end{align}
Now returning to \eqref{9.34}, invoking \eqref{eq7.9} and making note of \eqref{9.38} and \eqref{9.39}, 
we obtain, at the maximum point $(x_1, t_1)$, the inequality 
\begin{align}\label{7.32}
0 \ge &~F \Delta_\phi \zeta -2 (|\nabla \zeta |^2/\zeta) F + \zeta [\Delta_\phi -\partial_t] F \nonumber \\
 \ge & -F ([c_2 +(m-1) c_1(1+R \sqrt{k})]/R^2) -2c_1^2F/R^2 \nonumber\\
&+ \zeta [2t_1 (\Delta_\phi f)^2/m -F/t_1 -
2 \langle \nabla f, \nabla F \rangle - 2 t_1 (m-1)k |\nabla f|^2\nonumber\\
&+2 t_1 (\alpha -1)\langle \nabla f, \nabla (e^{-f}\mathscr G) \rangle +\alpha t_1 \Delta_\phi (e^{-f}\mathscr G)].
\end{align}

We point out that here and below we abbreviate the arguments of $\mathscr G=\mathscr G(t,x,e^f)$ for convenience. 
Next, in view of $\nabla (\zeta F) =0$ at $(x_1, t_1)$, that 
$\zeta \langle \nabla f, \nabla F \rangle = -F \langle \nabla f , \nabla \zeta \rangle$. 
Substituting these in \eqref{7.32} we can write after rearranging terms
\begin{align}\label{7.34}
0 \ge & -F [c_2 +(m-1) c_1(1+R \sqrt{k})+ 2c_1^2]/R^2 \nonumber\\
&-\zeta F/t_1 + 2 t_1 \zeta/m \left[|\nabla f|^2 + e^{-f} \mathscr G - \partial_t f \right]^2 \nonumber\\
& + 2F \langle \nabla f, \nabla \zeta \rangle - 2t_1 \zeta (m-1) k |\nabla f|^2  \nonumber\\
&+2 t_1 \zeta (\alpha -1)\langle \nabla f, \nabla (e^{-f}\mathscr G) \rangle +\alpha t_1 \zeta \Delta_\phi (e^{-f}\mathscr G).
\end{align} 
Now upon utilising the bound 
$2F \langle \nabla f, \nabla \zeta \rangle \ge -2 F |\nabla f| |\nabla \zeta| \ge -2F |\nabla f| (c_1/R) \sqrt \zeta$ 
it follows after multiplying \eqref{7.34} through by $t_1 \zeta (x_1)=t_1 \zeta$ and rearranging terms, that
\begin{align}\label{7.36}
0 \ge & -t_1 \zeta F ([c_2 +(m-1) c_1(1+R \sqrt{k})+ 2c_1^2]/R^2+ \zeta/t_1)\nonumber\\
& + 2 t_1^2 \zeta^2/m \left[|\nabla f|^2 + e^{-f} \mathscr G - \partial_t f \right]^2 \nonumber\\
& - 2t_1(c_1/R) \zeta^{3/2} |\nabla f| F- 2t_1^2 \zeta^2 (m-1) k |\nabla f|^2 \nonumber\\ 
&+t_1^2 \zeta^2 \left[2(\alpha -1)\langle \nabla f, \nabla (e^{-f}\mathscr G) \rangle+\alpha \Delta_\phi(e^{-f}\mathscr G)\right].
\end{align}

Referring next to the expression on the last line in \eqref{7.36} by working on the sum inside the brackets we can write  
\begin{align}\label{9.25}
2 (\alpha -1)\langle \nabla f,\nabla (e^{-f} \mathscr G) \rangle & + \alpha \Delta_\phi (e^{-f} \mathscr G)  \nonumber\\
=&~2 (\alpha -1) [ - e ^{-f} \mathscr G |\nabla f|^2 + e^{-f} \langle \nabla f, \nabla \mathscr G \rangle] \nonumber\\
& +\alpha [ e^{-f} \Delta_\phi  \mathscr G -2 e^{-f} \langle \nabla f , \nabla \mathscr G \rangle + \mathscr G \Delta_\phi e^{-f}] \nonumber\\
=&~-2 (\alpha -1) e ^{-f} \mathscr G |\nabla f|^2 + 2 \alpha e^{-f} \langle \nabla f , \nabla \mathscr G \rangle\nonumber\\
&-2e^{-f} \langle \nabla f, ( \mathscr G_x + e^ f \mathscr G_w \nabla f) \rangle 
+ \alpha e^{-f} [ \mathscr G_{xx} -\langle \nabla \phi, \mathscr G_x \rangle \nonumber\\
&+2e^f \langle \mathscr G_{xw} ,\nabla f \rangle 
+ e^f |\nabla f|^2 ( \mathscr G_w +e^f \mathscr G_{ww})+e^f \mathscr G_w \Delta_\phi f] \nonumber\\
& + \alpha \mathscr G e^{-f} ( -\Delta_\phi f+ |\nabla f|^2) -2 \alpha e^{-f} \langle \nabla f , \nabla \mathscr G \rangle.
\end{align}
Here we have made note of the relation 
\begin {align} \label{3.36}
\Delta_\phi e^{-f} &= \Delta e^{-f} - \langle \nabla \phi , \nabla e^{-f} \rangle =
{\rm div}(-e^{-f}\nabla f)+e^{-f} \langle \nabla \phi, \nabla f \rangle \nonumber\\
&= -e^{-f} \Delta f+ e^{-f} |\nabla f|^2 +e^{-f} \langle \nabla \phi , \nabla f\rangle 
= -e^{-f} (\Delta_\phi f- |\nabla f|^2).
\end{align}  
Next, since according to \eqref{9.6a} we have, 
\begin{align}
\Delta_\phi f [ \alpha \mathscr G_w - \alpha \mathscr G e^{-f}]&= 
\left[-\frac{F}{\alpha t_1}-\frac{\alpha -1}{\alpha}|\nabla f|^2 \right]
\left[ \alpha (\mathscr G_w - \mathscr G e^{-f})\right] \nonumber\\
& = -F (\mathscr G_w -\mathscr G e^{-f})/t_1 -(\alpha -1)|\nabla f|^2(\mathscr G_w -\mathscr G e^{-f}), 
\end{align}  
upon substituting this back in \eqref{9.25} it follows that  
\begin{align*}
2 (\alpha -1)&\langle \nabla f,\nabla (e^{-f} \mathscr G) \rangle +\alpha \Delta_\phi (e^{-f} \mathscr G) \nonumber \\
=& -2 (\alpha -1) e^{-f} \mathscr G |\nabla f|^2 
-2 e^{-f}\langle \nabla f, \mathscr G_x \rangle  -2 \mathscr G_w |\nabla f|^2 \\
& + \alpha e^{-f} ( \mathscr G_{xx} - \langle \nabla \phi , \mathscr G_x \rangle) + 2 \alpha \langle \mathscr G_{xw} ,\nabla f \rangle 
+ \alpha |\nabla f|^2 \mathscr G_w +\alpha |\nabla f|^2 e^f\mathscr G_{ww} \\
&+ \alpha \mathscr G e^{-f} |\nabla f|^2 - F (\mathscr G_w -\mathscr G e^{-f})/t_1-(\alpha -1)|\nabla f|^2(\mathscr G_w +\mathscr G e^{-f}) \\
=&~|\nabla f|^2 [-2 (\alpha -1) e^{-f} \mathscr G -2 \mathscr G_w + \alpha \mathscr G_w + 
\alpha e^f \mathscr G_{ww} + \alpha e^{-f}\mathscr G -(\alpha -1) \mathscr G_w \\
&+ (\alpha -1) \mathscr G e^{-f}] - [2 \langle \nabla f , (e^{-f}\mathscr G_x -\alpha \mathscr G_{xw})\rangle
 + F (\mathscr G_w -\mathscr G e^{-f})/t_1 -\alpha e^{-f}\Delta_\phi \mathscr G^x]. 
\end{align*}

Therefore, by taking into account the relevant cancellations, after simplifying terms and using basic inequalities, we can write  
\begin{align}\label{3.39}
2 (\alpha -1)&\langle \nabla f, \nabla (e^{-f}\mathscr G ) \rangle +\alpha \Delta_\phi(e^{-f}\mathscr G)\nonumber\\
\ge &~ |\nabla f|^2(e^{-f}\mathscr G - \mathscr G_w + \alpha e^f \mathscr G_{ww})  -F (\mathscr G_w -\mathscr G e^{-f})/t_1\nonumber\\
&-2 |\nabla f| |e^{-f} \mathscr G_x - \alpha \mathscr G_{xw}|+\alpha e^{-f} \Delta_\phi \mathscr G^x.
\end{align}

As a result making use of the relations \eqref{3.31}-\eqref{3.33} and the inequality \eqref{3.39} above and substituting 
all back into \eqref {7.36} and recalling $0 \le \zeta \le 1$ we obtain:
 \begin{align}\label{7.40}
0 \ge& - \zeta F ([c_2 +(m-1) c_1(1+R \sqrt{k})+ 2c_1^2]t_1 /R^2+ 1) 
- t_1\zeta^2 F \left[ \mathscr G_w - \mathscr G e^{-f} \right]\nonumber\\
&+ 2 t_1^2 \zeta^2/m \left[|\nabla f|^2 + e^{-f} \mathscr G - \partial_t f \right]^2 -2 c_1 t_1 \zeta^{3/2} |\nabla f| F/R \nonumber\\
&+ t_1^2 \zeta^2 |\nabla f|^2 \left[e^{-f}\mathscr G - \mathscr G_w + \alpha e^f \mathscr G_{ww} - 2(m-1)k \right]\nonumber\\
&-  2t_1^2 \zeta^2 |\nabla f| |e^{-f} \mathscr G_x - \alpha \mathscr G_{xw}| +\alpha t_1^2 \zeta^2 e^{-f}\Delta_\phi\mathscr G ^x.
\end{align}

In order to obtain the desired bounds out of this it is more efficient to introduce the quantities 
$y, z$ and $\mathsf{a}, \mathsf{b}$ by setting 
\begin{align} \label{sum-3.38}
\begin{cases}
y= \zeta |\nabla f|^2, \\ 
z= \zeta (\partial_t f -e^{-f}\mathscr G), \\
y-\alpha z= \zeta F/t_1>0, \\
\mathsf{a} = (m-1)k + (1/2) \sup \{[-e^{-f}\mathscr G + \mathscr G_w - \alpha e^f \mathscr G_{ww}]_+ 
: (x,t) \in H_{2R, T}\}, \\ 
\mathsf{b} = \sup \{|e^{-f} \mathscr G_x - \alpha \mathscr G_{xw}| : (x,t) \in H_{2R, T}\}. 
\end{cases}
\end{align}

Now substituting the quantities \eqref{sum-3.38} back in \eqref{7.40} and recalling $0 \le \zeta \le 1$, it follows that 
\begin{align} \label{7.43}
0 \ge& -\zeta F ([c_2 +(m-1) c_1(1+R \sqrt{k})+ 2c_1^2]t_1/R^2 + 1) \nonumber\\
&+ 2 t_1^2 /m \left[ (y-z)^2 -m c_1/R \sqrt{y}(y-\alpha z) -m\mathsf{a} y -m\mathsf{b} \sqrt y\right]\nonumber\\
& - t_1\zeta F \left[ \mathscr G_w - e^{-f}\mathscr G \right]_{+} +\alpha t_1^2 \zeta^2 [e^{-f}\Delta_\phi \mathscr G ^x]_-.
\end{align}   
In order to proceed further we now state and prove the following lemma.

\begin{lemma} \label{LiYau-basiclemma}
Suppose $z \in \mathbb{R}$, $c, y>0$ and $\alpha>1$ are arbitrary constants such that $y-\alpha z >0$. Then for any 
$\varepsilon \in (0,1)$ we have
\begin{align} \label{eq6.9} 
(y-z)^2& - (mc_1/R) \sqrt y (y-\alpha z) - m\mathsf{a} y - m\mathsf{b} \sqrt y \nonumber\\
&\ge (y-\alpha z)^2/\alpha^2-(mc_1/R)^2\alpha^2 (y-\alpha z)/[8 (\alpha-1)] \nonumber \\
& - (\alpha^2 m^2\mathsf{a}^2)/[4(1-\varepsilon)(\alpha-1)^2]
-(3/4) [m^4\mathsf{b}^4\alpha^2/(4 \varepsilon (\alpha -1)^2)]^{1/3}.  
\end{align}
\end{lemma}

\begin{proof}
Starting from the expression on the left-hand side in \eqref{eq6.9} we can write for any $\delta, \varepsilon$ by basic considerations
\begin {align} \label{eq6.10}
(y-z)^2 &- (mc_1/R) \sqrt y (y-\alpha z) - m\mathsf{a} y - m\mathsf{b} \sqrt y\nonumber\\
=&~(1-\varepsilon-\delta)y^2-(2-\varepsilon \alpha)yz+z^2+(\varepsilon y - (mc_1/R) \sqrt y)(y-\alpha z)\nonumber\\ 
& + \delta y^2 - m\mathsf{a}y - m\mathsf{b} \sqrt y \nonumber\\
=&~(1/\alpha - \varepsilon/2)(y-\alpha z)^2+(1-\varepsilon-\delta-1/\alpha+\varepsilon/2)y^2
+(1- \alpha + \varepsilon \alpha^2/2)z^2\nonumber\\
&+ (\varepsilon y - (mc_1/R) \sqrt y)(y-\alpha z) + \delta y^2 - m\mathsf{a}y - m\mathsf{b} \sqrt y. 
\end{align}   
In particular setting $\delta = (1/\alpha-1)^2$ and $\varepsilon = 2-2/\alpha-2(1/\alpha-1)^2 = 2(\alpha-1)/\alpha^2$ 
gives $1-\varepsilon-\delta-1/\alpha+\varepsilon/2=0$ and $1-\alpha + \varepsilon\alpha^2/2=0$ and so by making note 
of the inequality $\varepsilon y - (mc_1/R) \sqrt y \ge - (mc_1/R)^2/(4 \varepsilon)$ with $\varepsilon=2(\alpha-1)/\alpha^2>0$ we can 
deduce from 
\eqref{eq6.10} that
\begin {align}\label{8.36}   
(y-z)^2 - (mc_1/R) \sqrt y &(y-\alpha z) - m\mathsf{a} y - m\mathsf{b} \sqrt y \\
\ge&~ (y-\alpha z)^2/\alpha^2- m^2c_1^2\alpha^2 (y-\alpha z)/[8R^2 (\alpha-1)] \nonumber\\
&+ (\alpha -1)^2 y^2/\alpha^2 - m\mathsf{a}y- m\mathsf{b} \sqrt y. \nonumber 
\end{align} 
Next, considering the last three terms only we can write, for any $\varepsilon \in (0,1)$,
\begin {align} 
(\alpha -1)^2 &y^2/\alpha^2 - m\mathsf{a}y - m\mathsf{b} \sqrt y \nonumber\\
\ge&~ (\alpha -1)^2y^2/\alpha^2 - (1-\varepsilon)(\alpha-1)^2 y^2/\alpha^2 \nonumber\\
&-(\alpha^2 m^2\mathsf{a}^2)/[4(1-\varepsilon)(\alpha-1)^2] - m\mathsf{b} \sqrt y\nonumber\\
\ge&~ \varepsilon (\alpha -1)^2y^2/\alpha^2 - 
(\alpha^2 m^2\mathsf{a}^2)/[4(1-\varepsilon)(\alpha-1)^2] - m\mathsf{b} \sqrt y\nonumber\\
\ge& - (\alpha^2 m^2\mathsf{a}^2)/[4(1-\varepsilon)(\alpha-1)^2] 
 -(3/4) [m^4\mathsf{b}^4\alpha^2/(4 \varepsilon (\alpha -1)^2)]^{1/3} 
\end{align} 
where above we have made use of $(1-\varepsilon)(\alpha-1)^2 y^2/\alpha^2 - by \ge -(\alpha^2 b^2)/[4(1-\varepsilon)(\alpha-1)^2]$ 
and $\varepsilon (\alpha -1)^2y^2/\alpha^2 - c \sqrt y \ge  -(3/4) c^{4/3} [\alpha^2/(4 \varepsilon (\alpha -1)^2)]^{1/3}$  
to deduce the first and last inequalities respectively. Substituting back in \eqref{8.36} gives the desired inequality.  
\end{proof}

Now returning to the inequality \eqref{7.43} and making use of Lemma \ref{LiYau-basiclemma} it follows that   
\begin{align}\label{7.54}
0 \ge&-\zeta F ( [c_2 +(m-1) c_1(1+R \sqrt{k})+ 2c_1^2]t_1/R^2 + 1) \nonumber\\
&+ 2 t_1^2 /m [ (\zeta F)^2/t_1 ^2 \alpha^2 -m^2c_1^2 \alpha^2(\zeta F)/(8(\alpha-1)R^2 t_1)] \nonumber\\
&- m t_1 ^2\alpha^2\mathsf{a}^2/[2(1-\epsilon)(\alpha-1)^2]-[3t_1^2/2](m \alpha^2 \mathsf{b}^4/[4 \epsilon(\alpha-1)^2])^{1/3}\nonumber\\
&- t_1 \zeta F \left[ \mathscr G_w - e^{-f}\mathscr G \right]_{+}  +\alpha t_1^2 \zeta^2 [e^{-f}\Delta_\phi \mathscr G ^x]_-.
\end{align}

We now aim to rewrite \eqref{7.54} as an inequality involving a quadratic expression in powers of the localised function $\zeta F$. 
To this end upon setting
\begin{align} \label{7.55}
\mathsf{d} =& [c_2 +(m-1) c_1(1+R \sqrt{k})+2c_1^2]t_1/R^2 +1 \nonumber \\
&+ m t_1 c_1 ^2 \alpha ^2/[4(\alpha -1)R^2] + t_1\gamma^\mathscr G_{\mathsf C}(2R), 
\end{align}
and 
\begin{align} \label{Eeq6.61}
\mathsf{e}=&~m \alpha^2 \mathsf{a}^2/[2(1-\varepsilon)(\alpha-1)^2] \nonumber \\
&+ (3/2) [ m \alpha^2 \mathsf{b}^4/(4 \varepsilon(\alpha-1)^2) ]^{1/3} + \alpha \gamma^\mathscr G_{\mathsf D}(2R), 
\end{align}
where 
\begin{equation}
\gamma^\mathscr G_{\mathsf C} (2R) = \sup_{H_{2R, T}} \{[\mathscr G_w - e^{-f}\mathscr G]_{+}\}, \qquad 
\gamma^\mathscr G_{\mathsf D} (2R) = \sup_{H_{2R, T}} \{[-e^{-f}\Delta_\phi \mathscr G^x]_+\}, 
\end{equation}
it is seen that \eqref{7.54} can be expressed as 
\begin{align} \label{eq6.62}
0 \ge 2(\zeta F)^2/(m\alpha^2) - (\zeta F) \mathsf{d} - t_1^2 \mathsf{e}.
\end{align}
Now basic considerations using \eqref{eq6.62} on quadratics result in 
\begin {align} \label{eq6.63}
\zeta F &\le (m \alpha^2/4) \left(  \mathsf{d} + \sqrt{\mathsf{d}^2+ 8 t^2_1 \mathsf{e}/(m\alpha^2)} \right) \nonumber\\
&\le (m \alpha^2/4) \left( 2\mathsf{d}+\sqrt{(8t_1^2 \mathsf{e})/(m \alpha^2)} \right) = m \alpha^2\mathsf{d}/2 +t_1 \alpha \sqrt{m\mathsf{e}/2}.
\end{align}
Since $\zeta \equiv 1$ for $d(x,x_0) \le R$ and $(x_1, t_1)$ is the point where $\zeta F$ attains its maximum 
on $d(x,x_0) \le 2R$ we have
\begin{align}
F(x, \tau) = [\zeta F] (x, \tau) \le [\zeta F] (x_1,t_1) \le m \alpha^2\mathsf{d}/2 +t_1 \alpha \sqrt{m\mathsf{e}/2}.
\end{align}
Therefore recalling \eqref{2.3}, substituting for $\mathsf{d}$ and $\mathsf{e}$ from \eqref{7.55} and \eqref{Eeq6.61} above and 
making noting $t_1 \le \tau$, we can write after dividing both sides $\alpha \tau$, 
\begin{align}
\alpha^{-1} |\nabla f|^2 - \partial_t f + e^{-f} \mathscr G
\le&~(m \alpha/2\tau) \mathsf{d} + \sqrt{m \mathsf{e}/2}  \nonumber \\
\le&~(m \alpha) [c_2 +(m-1) c_1(1+R \sqrt{k})+2c_1^2] /2R^2 \nonumber \\
&+ (m \alpha/2\tau) 
+ (m \alpha/2) (\gamma^\mathscr G_{\mathsf C}(2R) + m c_1 ^2 \alpha ^2/[4(\alpha -1)R^2]) \nonumber \\
& + \sqrt {m/2} \{ m \alpha^2 \mathsf{a}^2/[2(1-\varepsilon)(\alpha-1)^2] \nonumber \\
& + (3/2) [ m \alpha^2 \mathsf{b}^4/(4 \varepsilon(\alpha-1)^2) ]^{1/3} + \alpha \gamma^\mathscr G_{\mathsf D}(2R) \}^{1/2}. 
\end{align}
Finally using the arbitrariness of $0< \tau \le T$ it follows after reverting back to $w$ upon noting the relation $f=\log w$ 
and rearranging terms that 
\begin{align} \label{Eq-3.46-Final}
\frac{|\nabla w|^2}{\alpha w^2} - \frac{\partial_t w}{w} + \frac{\mathscr G}{w} 
\le&~(m \alpha/2) [1/t + \gamma^\mathscr G_{\mathsf C}(2R)] \nonumber \\
&+ (m \alpha/2) [m c_1 ^2 \alpha ^2/[4(\alpha -1)]+c_2 \nonumber\\
&+(m-1) c_1(1+R \sqrt{k})+2c_1^2]/R^2 \nonumber \\
&+ \sqrt {m/2} \{ m \alpha^2 \mathsf{a}^2/[2(1-\epsilon)(\alpha-1)^2]\nonumber\\
&+ (3/2) [ m \alpha^2 \mathsf{b}^4/(4 \varepsilon(\alpha-1)^2) ]^{1/3}+ \alpha \gamma^\mathscr G_{\mathsf D}(2R) \}^{1/2}. 
\end{align}
Using \eqref{Eq-1.5}-\eqref{Eq-1.6} in conjunction with the bounds \eqref{eq2.2}-\eqref{eq2.3} and making note 
of \eqref{sum-3.38} it is seen that 
\begin{equation}
\mathsf{a} = (m-1) k + \gamma^{{\mathscr G}, \alpha}_\mathsf{A} (2R)/2, \qquad 
\mathsf{b} = \gamma^{{\mathscr G}, \alpha}_{\mathsf B}(2R)
\end{equation}
which upon substitution in \eqref{Eq-3.46-Final} gives the desired estimate as formulated in \eqref{1.26}. \hfill $\square$

\section{Proof of the parabolic Harnack inequality in Theorem \ref{thm38}} \label{sec8}

This will be shown as a consequence of the local Li-Yau estimate in Theorem \ref{thm7.1}. Here, the main idea is to integrate 
\eqref{1.26} along suitable space-times curves described below. Towards this end and upon referring 
to the expression on the right-hand side of \eqref{1.26} let us introduce the constant  
\begin{align} \label{H-in-Harnack-Eq}
\mathsf{H} =&~ \gamma^\mathscr G_\mathsf E (R) -   \frac{m\alpha}{2R^2} \left[ \frac{m c_1 ^2 \alpha ^2}{4(\alpha -1)} 
+c_2+(m-1) c_1(1+R \sqrt{k})+2c_1^2 \right]\nonumber\\
&- (m \alpha/2) \gamma^{\mathscr G}_{\mathsf C} (2R) 
- \sqrt {\frac{m \alpha}{2}}\bigg\{ \frac{m \alpha [(m-1) k 
+ \gamma^{{\mathscr G}, \alpha}_\mathsf{A} (2R)/2]^2}{2(1-\varepsilon)(\alpha-1)^2}\nonumber\\
&+ \left[\frac{3^3 m [\gamma^{{\mathscr G}, \alpha}_\mathsf{B} (2R)]^{4}}{2^5 \varepsilon \alpha (\alpha-1)^2}\right]^{1/3} 
+ \gamma^{\mathscr G}_{\mathsf D} (2R) \bigg\}^{1/2},
\end{align}
where in analogy with the other $\gamma$-quantities occurring earlier, the first term on the right stands for 
\begin{equation} \label{gamma-in-Harnack-Eq}
\gamma^\mathscr G_\mathsf E(R) = \inf_{H_{R,T}} w^{-1} \mathscr G(t,x,w).
\end{equation}
It then follows from Theorem \ref{thm7.1} that in the space-time cylinder $H_{R,T}$ we have the inequality  
\begin{align} \label{Eq-4.3}
\frac{\partial_t w}{w} & \ge \frac{|\nabla w|^2}{\alpha w^2} - \frac{m \alpha}{2t} + \mathsf{H}. 
\end{align}

Note that in \eqref{H-in-Harnack-Eq} we can replace $\gamma^\mathscr G_\mathsf E(R)$ 
with the smaller quantity $\gamma^\mathscr G_\mathsf E(2R)$ in order to maintain uniformity 
in notation (specifically, have all the $\gamma$-quantities in the formula defining $\mathsf H$ with argument $2R$). However, as far as the bound 
\eqref{Eq-4.3} is concerned, the larger quantity $\gamma^\mathscr G_\mathsf E(R)$ is sufficient and clearly more accurate. 

Next suppose $\zeta \in \mathscr{C}^1( [t_1,t_2]; M)$ is an arbitrary curve in $M$ lying entirely in ${\mathscr B}_R$ and 
satisfying $\zeta(t_1) = x_1$ and $\zeta(t_2) = x_2$. In particular $(\zeta(t), t) \in H_{R,T}$ for all $t_1 \le t \le t_2$. 
Differentiating the function $\log w(\zeta(t), t)$ and using \eqref{Eq-4.3} it is then seen that 
\begin{align}
d/dt [\log w(\zeta(t),t)] &= \langle \nabla w/w, \dot \zeta (t) \rangle + \partial_t w/w \nonumber \\
&\ge  \langle \nabla w/w, \dot \zeta (t) \rangle + |\nabla w|^2/(\alpha w^2) - ( m \alpha)/(2t) + \mathsf{H} \nonumber\\
&= \alpha^{-1} |\nabla w/w + \alpha \dot \zeta (t)/2|^2 - \alpha |\dot \zeta (t)|^2/4 
-(m \alpha)/(2t) + \mathsf{H} \nonumber\\ 
&\ge -\alpha |\dot \zeta (t)|^2/4 - (m \alpha)/(2t) + \mathsf{H}. 
\end{align}
Therefore integrating the above inequality gives 
\begin{align}
\log \frac{w(x_2,t_2)}{w(x_1,t_1)} &= \log w(\zeta(t),t)\bigg|_{t_1}^{t_2}= \int_{t_1}^{t_2} \frac{d}{dt} \log w(\zeta(t),t) \, dt \nonumber \\
&\ge \int_{t_1}^{t_2}  -\frac{\alpha}{4} |\dot \zeta (t)|^2 \,dt 
- \int_{t_1}^{t_2} \frac{m \alpha}{2t} dt+\int_{t_1}^{t_2} \mathsf{H} \, dt \nonumber \\
&= - (\alpha/4) \int_{t_1}^{t_2} |\dot \zeta (t)|^2 \,dt  -m \alpha/2 \log (t_2/t_1) + (t_2-t_1) \mathsf{H}.
\end{align}
Hence upon exponentiating we have
\begin{align}
\frac{w(x_2,t_2)}{w(x_1,t_1)} &\ge  
\exp \left[ - \int_{t_1}^{t_2} \frac{\alpha}{4}|\dot \zeta(t)|^2 \, dt \right]  \exp [(t_2-t_1) \mathsf{H}] \left(\frac{t_2}{t_1} \right)^{-m \alpha/2} ,  
\end{align}
or after rearranging terms and rescaling the integral:  
\begin{align}
\frac{w(x_2,t_2)}{w(x_1,t_1)} \ge {\rm exp} [ (t_2-t_1) \mathsf{H}  -\alpha L(x_1,x_2, t_2-t_1)]  \left(\frac{t_2}{t_1}\right)^{-m \alpha/2}
\end{align}
where 
\begin{align}
L(x_1,x_2, t_2-t_1)=\inf_{\zeta \in \Gamma} \left[ \frac{1}{4(t_2-t_1)} \int_{0}^{1} |\dot \zeta(t)|^2\,dt \right].
\end{align}
This gives the parabolic Harnack inequality in its local form. Now if the bounds are global as in Theorem \ref{thm28-global} 
then by passing to the limit $R \to \infty$ we obtain the global counterpart of the inequality with the constant ${\mathsf H}$ 
adjusted. \hfill $\square$

\section{Proof of the Liouville Theorems and implications $({\bf II})$: ${\mathscr Ric}_\phi^m(g) \ge 0$} \label{sec9}

In this section we give the proofs of the Liouville-type results formulated in Theorem \ref{coroLiouville} and its implications 
as given in Theorems \ref{LiouvilleThmEx} and \ref{LiouvilleThmEx-log}. 
The basis for all these proofs is the elliptic estimate in Theorem \ref{thm48} below 
that in term is a consequence of the parabolic Li-Yau estimates in Theorems \ref{thm7.1} and \ref{thm28-global}. 
Towards this end we begin by presenting the proof of the estimate \eqref{eqL2.26}. 
We emphasise that in all these results we are assuming the curvature condition ${\mathscr Ric}_\phi^m(g) \ge 0$ 
with $n \le m <\infty$.

\qquad \\
{\bf Proof of Theorem  \ref{thm48}} The proof of \eqref{eqL2.26} for the elliptic equation $\Delta_\phi w + {\mathscr G}(w) =0$ follows 
directly from the global estimate \eqref{1.26-global} 
in Theorem \ref{thm28-global} by passing to the limit $t \nearrow \infty$, upon noting that $w$ is time independent. 
Moreover, here, ${\mathscr G}$ is independent of $x$ and so by referring to \eqref{Eq-1.6} and \eqref{Eq-1.8} 
it follows that 
\begin{equation}
\mathsf{B}_{\mathscr G}^\alpha (w) \equiv 0, \qquad \mathsf{D}_{\mathscr G} (w) \equiv 0,  
\end{equation}
thus giving $\gamma^{{\mathscr G}, \alpha}_\mathsf{B} =0$ and $ \gamma^{\mathscr G}_\mathsf{D}=0$. As a result the middle 
expression ``$\mathsf M$" on the right in \eqref{1.26-global} reduces to 
\begin{align}
\mathsf{M} &= \sqrt {\frac{m \alpha}{2}}\bigg\{ \frac{m \alpha [(m-1) k 
+ \gamma^{{\mathscr G}, \alpha}_\mathsf{A}/2]^2}{2(1-\varepsilon)(\alpha-1)^2}
+ \left[\frac{3^3 m [\gamma^{{\mathscr G}, \alpha}_\mathsf{B}]^{4}}{2^5 \varepsilon \alpha(\alpha-1)^2}\right]^{1/3} 
+ \gamma^{\mathscr G}_{\mathsf D} \bigg\}^{1/2} \nonumber \\
&= \sqrt {\frac{m \alpha}{2}}\bigg\{ \frac{m\alpha [(m-1) k + \gamma^{{\mathscr G}, \alpha}_\mathsf{A}/2]^2}{2(1-\varepsilon)(\alpha-1)^2} \bigg\}^{1/2} 
= \frac{m \alpha}{4} \frac{[2(m-1) k + \gamma^{{\mathscr G}, \alpha}_\mathsf{A}]}{(\alpha-1) \sqrt{1-\varepsilon}}. 
\end{align}
Putting the above fragments together gives the desired estimate. \hfill $\square$

\qquad \\
{\bf Proof of Theorem \ref{coroLiouville}} Since here ${\mathscr Ric}_\phi^m (g) \ge 0$, by taking $k=0$ in \eqref{eqL2.26} 
it follows that under the assumptions of the theorem $w$ satisfies the global estimate on $M$:
\begin{align} \label{Eq-5.3}
\frac{|\nabla w|^2}{\alpha w^2} + \frac{\mathscr G{}(w)}{w} 
\le&~\frac{m \alpha}{4(\alpha-1)\sqrt{1-\varepsilon}} \gamma^{{\mathscr G}, \alpha}_\mathsf{A}
+ \frac{m \alpha}{2}  \gamma^{\mathscr G}_\mathsf{C}.
\end{align}
Next, since by assumption ${\mathscr G}(w) - w {\mathscr G}_w(w) \ge 0$ and 
${\mathscr G}(w)-w{\mathscr G}_w(w)+\alpha w^2 {\mathscr G}_{ww}(w) \ge 0$, 
for some $\alpha>1$, it follows upon recalling \eqref{Eq-1.5} and \eqref{Eq-1.7} that 
\begin{align}
\mathsf{A}_{\mathscr G}^\alpha (w) 
= [- \alpha w{\mathscr G}_{ww}(w) + {\mathscr G}_w(w) - {w^{-1} \mathscr G}(w)]_+
\equiv 0, \nonumber \\
\mathsf{C}_{\mathscr G}^\alpha (w) 
= [{\mathscr G}_w(w) - w^{-1}{\mathscr G}(w)]_{+}
\equiv 0, 
\end{align}
thus giving $\gamma^{{\mathscr G}, \alpha}_\mathsf{A} =0$ and $ \gamma^{\mathscr G}_\mathsf{C}=0$. 
Substituting these back in \eqref{Eq-5.3} then leads to 
\begin{align} 
\frac{|\nabla w|^2}{\alpha w^2} + \frac{\mathscr G(w)}{w} \le 0.  
\end{align}
Since $\mathscr G(w) \ge 0$ we thus infer that $|\nabla w|^2/(\alpha w^2) + \mathscr G(w)/w \equiv 0$ 
and therefore $|\nabla w| \equiv 0$ on $M$. The conclusion on $w$ being a constant is now immediate. 
\hfill $\square$

\qquad \\
{\bf Proof of Theorem \ref{LiouvilleThmEx}} 
A direct calculation gives  
\begin{align}
{\mathscr G} (w) - w {\mathscr G}_w (w) &= \sum_{j=1}^N \mathsf{A}_j (1-p_j ) w^{p_j}, \\
{\mathscr G} (w) - w {\mathscr G}_w (w) + \alpha w^2 {\mathscr G}_{ww} (w) 
&= \sum_{j=1}^N [\alpha \mathsf{A}_j p_j (p_j -1) - \mathsf{A}_j p_j + \mathsf{A}_j] w^{p_j} \nonumber \\
&= \sum_{j=1}^N [\mathsf{A}_j (p_j -1) (\alpha p_j -1)] w^{p_j}. 
\end{align}

Now if $\mathsf{A}_j \ge 0$ we have ${\mathscr G}(w) \ge 0$ and the assumption $p_j \le 1$ gives 
${\mathscr G} (w) - w {\mathscr G}_w (w) \ge 0$ and 
${\mathscr G}(w) - w {\mathscr G}_w (w) + \alpha w^2 {\mathscr G}_{ww} (w) \ge 0$ (by 
choosing $\alpha>1$ suitably). 
\hfill $\square$

\qquad \\
{\bf Proof of Theorem \ref{LiouvilleThmEx-log}} 
Here  
${\mathscr G} (w)-w {\mathscr G}_w (w) = \mathsf{B} (1-s) w^{s} - \mathsf{A} w Y'(\log w)$ 
and ${\mathscr G} (w)- w {\mathscr G}_w (w) + \alpha w^2 {\mathscr G}_{ww} (w) 
= {\mathsf A}[(\alpha-1) w Y' + \alpha w Y''] + \mathsf{B} (s -1) (\alpha s -1) w^{s}$. 
Therefore the conclusion follows from an application of Theorem \ref{coroLiouville}.
\hfill $\square$

\qquad \\\
{\bf Data Availability Statement.} 
All data in this research is provided in full in the results section. 
Additional data is in the public domain at locations cited in the reference section.

\qquad \\
{\bf Acknowledgements.} The authors gratefully acknowledge financial support from the Engineering and Physical 
Sciences Research Council (EPSRC) through the research grant EP/V027115/1. They also wish to thank the 
anonymous reviewers for a careful reading of the paper and useful comments.

\end{document}